\input amstex
\input amsppt.sty
\magnification=\magstep1
\hsize=30truecc
\vsize=22.2truecm
\baselineskip=16truept
\TagsOnRight
\nologo
\pageno=1
\topmatter
\def\N{\Bbb N}
\def\Z{\Bbb Z}
\def\Q{\Bbb Q}

\def\l{\left}
\def\r{\right}
\def\b{\bigg}

\def\({\b(}
\def\[{\b[}
\def\){\b)}
\def\]{\b]}

\def\t{\text}
\def\f{\frac}
\def\mo{\roman{mod}}
\def\ord{\roman{ord}}

\def\bi{\binom}
\def\eq{\equiv}

\def\ls{\leqslant}
\def\gs{\geqslant}
\def\al{\alpha}

\def\Proof{\noindent{\it Proof}}
\def\Remark{\medskip\noindent{\it Remark}}
\def\Ack{\noindent {\bf Acknowledgment}}
\hbox {Preprint (August 14, 2006)}
\bigskip
\title Fleck quotients and Bernoulli numbers\endtitle
\author  Zhi-Wei Sun\endauthor
\affil Department of Mathematics, Nanjing University
\\ Nanjing 210093, People's Republic of China
\\zwsun\@nju.edu.cn
\\ {\tt http://pweb.nju.edu.cn/zwsun}
\medskip
\endaffil
\thanks 2000 {\it Mathematics Subject Classification}.\,Primary 11B65;
Secondary 05A10, 11A07, 11B68, 11B73, 11S80, 11T24.\newline\indent
Supported by the National Science Fund
for Distinguished Young Scholars (No. 10425103) in China.
\endthanks
\abstract Let $p$ be a prime, and let $n>0$ and $r$ be integers.
In 1913 Fleck showed that
$$F_p(n,r)=(-p)^{-\lfloor(n-1)/(p-1)\rfloor}\sum_{k\eq r\,(\mo\ p)}\bi nk(-1)^k\in\Z.$$
Nowadays this result plays important roles in many aspects.
Recently Sun and Wan investigated $F_p(n,r)$ mod $p$
in [SW2]. In this paper, using $p$-adic methods
we determine $(F_p(m,r)-F_p(n,r))/(m-n)$ modulo $p$
in terms of Bernoulli numbers, where $m>0$ is an integer
with $m\not=n$ and $m\eq n\ (\mo\ p(p-1))$.
Consequently, $F_p(n,r)$ mod $p^{\ord_p(n)+1}$ is determined;
for example, if $n\eq n_*\ (\mo\ p-1)$ with $0<n_*<p-2$
then
$$\f{F_p(pn,0)}{pn}\eq\f {n_*!}{n_*+1}B_{p-1-n_*}\ (\mo\ p).$$
This yields an application to Stirling numbers of the second kind.
We also study extended Fleck quotients; in particular we prove that
if $a>0$ and $l\gs0$ are integers with $2\ls n-l\ls p$ then
$$\f1{p^{n-l}}\sum_{l<k\ls n}\bi{p^an-d}{p^ak-d}(-1)^{pk}\bi{k-1}l
\eq\f{(-1)^{l-1}n!}{l!(n-l)}B_{p-n+l}\ (\mo\ p)$$
for all $d=1,\ldots,\max\{p^{a-2},1\}$.
\endabstract
\endtopmatter

\document

\heading{1. Introduction}\endheading

Let $p$ be a prime.
In 1819 C. Babbage observed that
$$\align &(-1)^{p-1}\bi{2p-1}{p-1}=\prod_{k=1}^{p-1}\l(1-\f{2p}k\r)
\\\eq&1-2p\sum_{k=1}^{p-1}\f1k=1-p\sum_{k=1}^{p-1}\l(\f1k+\f1{p-k}\r)\eq1\ (\mo\ p^2).
\endalign$$
In 1862 J. Wolstenholme proved further that if $p>3$ then
$$\bi{2p-1}{p-1}=\f12\bi{2p}p\eq1\ (\mo\ p^3).$$
This is a fundamental congruence involving binomial coefficients.
When $p>3$, we also have
$$(-1)^{r+1}\bi{2p-1}{p-r-1}\eq-2p^2H_r^2+2pH_r-1\ (\mo\ p^3)\tag1.1$$
for each $r=1,\ldots,p-1$, where $H_r=\sum_{0<k\ls r}1/k$;
in fact,
$$\align&(-1)^{r+1}\bi{2p-1}{p-r-1}+\bi{2p-1}{p-1}
\\=&\sum_{s=1}^r\((-1)^{s+1}\bi{2p-1}{p-s-1}-(-1)^s\bi{2p-1}{p-s}\)
\\=&\sum_{s=1}^r(-1)^{s+1}\f{2p}{p-s}\bi{2p-1}{p-s-1}
=2p\sum_{s=1}^r\f{s+p}{s^2-p^2}\prod_{0<k<p-s}\l(1-\f{2p}k\r)
\\\eq&2p\sum_{s=1}^r\l(\f1s+\f p{s^2}\r)\(1-2p\(H_{p-1}-\sum_{t=1}^s\f1{p-t}\)\)
\\\eq&2p\sum_{s=1}^r\l(\f1s+\f p{s^2}\r)(1-2pH_s)\eq2p\sum_{s=1}^r\l(\f1s+\f p{s^2}-2p\f{H_s}s\r)
\ (\mo\ p^3)
\endalign$$
and hence
$$\align(-1)^{r+1}\bi{2p-1}{p-r-1}+1\eq&2pH_r+2p^2\sum_{s=1}^r\f1{s^2}-4p^2\sum_{1\ls t\ls s\ls r}\f1{st}
\\\eq&2pH_r-2p^2H_r^2\ (\mo\ p^3).
\endalign$$

Let $n\in\N=\{0,1,2,\ldots\}$ and $r\in\Z$.
In 1913 A. Fleck (cf. [D, p.\,274]) showed that
$$\ord_p(C_p(n,r))\gs\l\lfloor\f{n-1}{p-1}\r\rfloor,$$
where $\lfloor\cdot\rfloor$
is the well-known floor function, $\ord_p(\al)$ denotes the $p$-adic order
of a $p$-adic number $\al$ (we regard $\ord_p(0)$ as $+\infty$), and
$$C_p(n,r)=\sum_{k\eq r\,(\mo\ p)}\bi nk(-1)^k.\tag1.2$$
(In [S02] the author expressed certain sums like (1.2) in terms of linear recurrences.)
Fleck's result plays fundamental roles in the recent investigation of the $\psi$-operator related
to Fontaine's theory, Iwasawa's theory and $p$-adic Langlands
correspondence (cf. [Co], [W] and [SW1]), and Davis and Sun's study of
homotopy exponents of special unitary groups (cf. [DS] and [SD]).
It is also related to Leopoldt's formula for $p$-adic $L$-functions (cf. [Mu, Theorem 8.5]).

 Note that if $p\not=2$ then
$$C_p(2p-1,-1)=\bi{2p-1}{p-1}(-1)^{p-1}+\bi{2p-1}{2p-1}(-1)^{2p-1}=\bi{2p-1}{p-1}-1$$
and $$C_p(2p,0)=\bi{2p}0(-1)^0+\bi{2p}p(-1)^p+\bi{2p}{2p}(-1)^{2p}=2-\bi{2p}p.$$
So, in the case $n=2p-1$ and $r=-1$,
or the case $n=2p$ and $r=0$, Fleck's result yields Babbage's congruence.

For $m=0,1,2,\ldots$, the $m$th order Bernoulli polynomials
$B^{(m)}_k(t)\  (k\in\N)$ are given by
$$\f{x^me^{tx}}{(e^x-1)^m}=\sum_{k=0}^\infty B^{(m)}_k(t)\f{x^k}{k!}\ \ (0<|x|<2\pi),$$
and those $B_k^{(m)}=B^{(m)}_k(0)\ (k\in\N)$ are called $m$th order Bernoulli numbers.
Clearly $B_k^{(0)}(t)=t^k$.
The usual Bernoulli polynomials and numbers are $B_k(t)=B_k^{(1)}(t)$
and $B_k=B_k(0)=B_k^{(1)}$ respectively.
It can be easily seen that
$$B^{(m)}_k(t)=\sum_{j=0}^k\bi kjB^{(m)}_jt^{k-j}\ \ \t{and}\ \ B_k^{(m)}(m-t)=(-1)^kB_k^{(m)}(t).$$
Since $B^{(m)}_k/k!$ coincides with $[x^k](x/(e^x-1))^m$, the coefficient
of $x^k$ in the power series expansion of $(x/(e^x-1))^m$, if $m\in\Z^+=\{1,2,3,\ldots\}$ then
$$\f{B^{(m)}_k}{k!}=\sum_{i_1+\cdots+i_m=k}
\f{B_{i_1}\cdots B_{i_m}}{i_1!\cdots i_m!}.$$
It is well known that $B_0=1$, $B_1=-1/2$ and $B_{2k+1}=0$ for $k=1,2,3,\ldots$.
The von Staudt-Clasusen theorem (cf. [IR, pp.\,233--236] or [Mu, Theorem 2.7])
states that $B_{2k}+\sum_{p-1\mid 2k}1/p\in\Z$
for any $k\in\Z^+$. Thus, $B_0,\ldots,B_{p-2}$ are $p$-adic integers
and hence so are $B_0^{(m)},\ldots,B_{p-2}^{(m)}$. Therefore
$B^{(m)}_k(t)\in\Z_p[t]$ if $0\ls k<p-1$, where $\Z_p$ denotes the ring of $p$-adic integers.

In terms of higher-order Bernoulli polynomials,
the author and D. Wan [SW2] determined the Fleck quotient
$$F_{p}(n,r):=(-p)^{-\lfloor(n-1)/(p-1)\rfloor}\sum_{k\eq r\,(\mo\ p)}\bi nk(-1)^k+[\![n=0]\!]\tag1.3$$
modulo $p$. (Throughout this paper, for an assertion
$A$ we let $[\![A]\!]$ take $1$ or 0 according as $A$ holds or not.)
Namely, by [SW2, Theorem 1.2], we have
$$F_p(n,r)\eq-n_*!B_{n^*}^{(m)}(-r)\ (\mo\ p)\ \ \t{for}\ m\in\N\ \t{with}\ m\eq-n\ (\mo\ p),\tag1.4$$
where $n_*$ is the smallest positive residue of $n$ modulo $p-1$,
and $n^*$ is the least nonnegative residue $\{-n\}_{p-1}$ of $-n$ modulo $p-1$.
For convenience the notations $n_*$ and $n^*$ will be often used, and we remind the reader
of the difference. Note that $n_*+n^*=p-1$ and hence
$$n_*!n^*!=\f{n_*!(p-1)!}{\prod_{0<k\ls n_*}(p-k)}\eq(-1)^{n_*-1}\eq(-1)^{n^*-1}\eq(-1)^{n-1}\ (\mo\ p).$$
We mention that $((p-1)/2)!$ mod $p$ is related to the class number of the quadratic field
$\Q(\sqrt{(-1)^{(p-1)/2}p})$ (cf. [Ch] and [M]).

Let $\zeta_p$ be a fixed primitive $p$th root of unity in the algebraic closure
of the $p$-adic filed $\Q_p$. It is easy to see that $\ord_p(1-\zeta_p)=1/(p-1)$.
The main trick in [SW2] is to determine $F_p(n,r)$ modulo $1-\zeta_p$.

Corollary 1.5(i) of Sun and Wan [SW2] states that
if $2\ls n\ls p$ then
$$\f1{p^n}\sum_{k=1}^n\bi{pn-1}{pk-1}(-1)^{pk}\eq-(n-1)!B_{p-n}\ (\mo\ p).\tag1.5$$
This is a further extension of Wolstenholme's congruence.
When $1<n<p-1$ and $2\mid n$, the right-hand side of the above congruence is zero since
$B_{p-n}=0$; inspired by (1.5) the author's student H. Pan used {\it Mathematica}
to find the conjecture
$$\f1{p^n}\sum_{k=1}^n\bi{pn-1}{pk-1}(-1)^k
\eq\f{n!n}{2(n+1)}pB_{p-1-n}\ (\mo\ p^2).\tag1.6$$

For $m,n\in\Z^+$, the Stirling number $S(m,n)$ of the second kind
is the number of ways to partition a set of cardinality $m$
into $n$ subsets. It is easy to show that
$$(-1)^{n-1}n!S(m\varphi(p^b),n)\eq C_p(n,0)=(-p)^{\lfloor(n-1)/(p-1)\rfloor}F_p(n,0)\ (\mo\ p^b)$$
for any $b=1,2,3,\ldots$ (cf. [GL, Lemma 5]),
where $\varphi$ is Euler's totient function. Thus, for sufficiently large $b>0$,
we have
$$\ord_p(n!S(m\varphi(p^b),n))=\l\lfloor\f{n-1}{p-1}\r\rfloor+\ord_p(F_p(n,0)).$$

In this paper we want to reveal further connections between
Fleck quotients and Bernoulli numbers, including
the determination of $F_p(pn,r)$ modulo $p^{\ord_p(n)+2}$
by which (1.6) holds when $1<n<p-1$ and $2\mid n$.
The method of [SW2] does not work for this purpose;
instead of $1-\zeta_p$ we define
$$\pi:=-\sum_{k=1}^{p-1}\f{(1-\zeta_p)^k}{k}\in\Z_p[\zeta_p].\tag1.7$$
It can be shown that $\pi^{p-1}/p\eq-1\ (\mo\ p)$
(see Section 2). The $p$-adic method in this paper
deals with congruences modulo powers of $\pi$, and it is so powerful that
we need not appeal to the Stickelberger congruence (cf. Theorems 11.2.1 and 11.2.10 of [BEW])
which is of advanced nature.

Sun and Wan [SW2, Corollary 1.7] proved the following periodical results:
$$F_p(n+p^b(p-1),r)\eq F_p(n,r)\ (\mo\ p^b)\quad\t{for}\ b=1,2,3,\ldots.$$
Thus, if $m\in\N$, $m\not=n$ and $m\eq n\ (\mo\ p(p-1))$,
then $(F_p(m,r)-F_p(n,r))/(m-n)\in\Z_p$. We determine this quotient
modulo $p$ in our first theorem.

\proclaim{Theorem 1.1} Let $p$ be a prime, and let
$m,n\in\N$, $m\not=n$ and $m\eq n\ (\mo\ p(p-1))$.
Then we have
$$\aligned&\f{F_p(m,r)-F_p(n,r)}{m-n}
\eq\f{(-1)^{n^*}}{n^*!}\sum_{1<k\ls n^*}\bi{n^*}k\f{B_k}kB^{(\{-n\}_p)}_{n^*-k}(-r)
\\\eq&(-1)^{n_*-1}n_*!
\sum_{1<k\ls n^*}\bi{n_*+k}{n_*}r^{n^*-k}\sum_{1<j\ls k}\bi
kj\f{B_j}jB^{(\{-n\}_p)}_{k-j} \ (\mo\ p).
\endaligned\tag1.8$$
\endproclaim

\proclaim{Corollary 1.1} Let $p\gs5$ be a prime,
and let $n>0$ be an integer with $n\not\eq0,-1\ (\mo\ p-1)$.
Then, there are at least $p-n^*+2\gs 4$ values of $r\in\{0,1,\ldots,p-1\}$
such that $\ord_p(F_p(m,r)-F_p(n,r))=\ord_p(m-n)$
for all $m\in\N$ with $m\eq n\ (\mo\ p(p-1))$.
\endproclaim
\Proof. Clearly $n^*=\{-n\}_{p-1}\gs2$ and $n_*=p-1-n^*<p-2$. Since the polynomial
$$P(x)=\sum_{k=2}^{n^*}\bi{n_*+k}{n_*}x^{n^*-k}\sum_{j=2}^k\bi kj\f{B_j}jB_{k-j}^{(\{-n\}_p)}\in\Z_p[x]$$
has degree at most $n^*-2$, and
$$\align[x^{n^*-2}]P(x)=&\bi{n_*+2}{n_*}\sum_{j=2}^2\bi{2}j\f{B_j}jB_{2-j}^{(\{-n\}_p)}
\\=&\bi{n_*+2}{n_*}\f{B_2}2B_0^{(\{-n\}_p)}=\f{(n_*+1)(n_*+2)}{24}\not\eq0\ (\mo\ p),
\endalign$$
there are at most $n^*-2$ values of $r\in\{0,1,\ldots,p-1\}$
satisfying $P(r)\eq0\ (\mo\ p)$. (Recall that a polynomial of degree $d\gs0$
over a field cannot have more than $d$ zeroes in the field.)
Combining this with Theorem 1.1 we obtain the desired result. \qed

\proclaim{Corollary 1.2} Let $p$ be a prime and let $r\in\Z$.
Suppose that $n\in\N$ and $n\eq0\ (\mo\ p-1)$. Then, for any $b=2,3,\ldots$ we have
$$F_p(n,r)\eq F_p(\{n\}_{\varphi(p^b)},r)\ (\mo\ p^b)\tag1.9$$
and
$$F_p(n+p-2,r)\eq F_p(\{n\}_{\varphi(p^b)}+p-2,r)\ (\mo\ p^b).tag1.10$$
Consequently,
$$\f{F_p(pn,r)+p[\![p\mid r]\!]-1}{pn}\eq0\ (\mo\ p).\tag1.11$$
\endproclaim
\Proof. (i) Let $b>1$ be an integer. Write $n=\varphi(p^b)q+\{n\}_{\varphi(p^b)}$
with $q\in\N$. As $n^*=\{-n\}_{p-1}=0$ and $(n+p-2)^*=[\![p\not=2]\!]$, if $q>0$ then
by Theorem 1.1 we have
$$\f{F_p(n,r)-F_p(\{n\}_{\varphi(p^b)},r)}{p^{b-1}(p-1)q}\eq0\ (\mo\ p)$$
and
$$\f{F_p(n+p-2,r)-F_p(\{n\}_{\varphi(p^b)}+p-2,r)}{p^{b-1}(p-1)q}\eq0\ (\mo\ p).$$
So (1.9) and (1.10) are valid.

(ii) Clearly $b=\ord_p(pn)+1\gs2$ and $\varphi(p^b)\mid pn$.
In view of (1.9),
$$F_p(pn,r)\eq F_p(0,r)=-pC_p(0,r)+1=1-p[\![p\mid r]\!]\ (\mo\ p^b)$$
and hence (1.11) holds. \qed

\proclaim{Corollary 1.3} Let $p$ be a prime and let $r\in\Z$.
Assume that $m,n\in\N$, $m\not=n$ and $m\eq n\ (\mo\ p-1)$.
Then
$$\aligned&\f{F_p(pm+p-1,r)-F_p(pn+p-1,r)}{p(m-n)}
\\\eq&\f{(-1)^{n^*}}{n^*!}\(\sum_{0<k<n^*}\f{B_k(-r)}kB_{n^*-k}(-r)-H_{n^*-1}B_{n^*}(-r)\)
\ (\mo\ p).\endaligned\tag1.12$$
\endproclaim
\Proof. In light of Theorem 1.1,
$$\f{F_p(pm+p-1,r)-F_p(pn+p-1,r)}{p(m-n)}\eq\f{(-1)^{n^*}}{n^*!}B\ (\mo\ p),$$
where
$$B=\sum_{1<k\ls n^*}\bi{n^*}k\f{B_k}kB_{n^*-k}(-r).$$
By a polynomial form of Miki's identity (cf. [PS, (2.3)]),
$$\align B+H_{n^*-1}B_{n^*}(-r)=&\f{n^*}2\sum_{0<k<n^*}\f{B_k(-r)B_{n^*-k}(-r)}{k(n^*-k)}
\\=&\f12\sum_{0<k<n^*}\l(\f1k+\f1{n^*-k}\r)B_k(-r)B_{n^*-k}(-r)
\\=&\sum_{0<k<n^*}\f{B_k(-r)}kB_{n^*-k}(-r).
\endalign$$
So we have the desired (1.12). \qed

\proclaim{Lemma 1.1} Let $p$ be an odd prime, and let $n$ be a positive integer
not divisible by $p-1$. Then $F_p(n,0)/n\in\Z_p$.
\endproclaim

With help of this lemma and Theorem 1.1, we can deduce the following theorem.

\proclaim{Theorem 1.2} Let $p$ be an odd prime, and let $n\in\Z^+$ with $2\mid n$ and $p-1\nmid n$.
Then
$$\f{2\sum_{k=1}^n\bi{pn-1}{pk-1}(-1)^k}{(-p)^{\lfloor(n-2)/(p-1)\rfloor}p^{n+1}n}
=\f{F_p(pn,0)}{pn}\eq\f{n_*!}{n_*+1}B_{p-1-n_*}\ (\mo\ p).\tag1.13$$
Given $b,m\in\Z^+$ with $b>2\lfloor(pn-1)/(p-1)\rfloor$, we have
$$\aligned\f2{pn}\cdot\f{(pn-1)!S(m\varphi(p^b),pn-1)}{(-p)^{\lfloor(pn-2)/(p-1)\rfloor}}
\eq&-\f{(pn-1)!S(m\varphi(p^b),pn)}{(-p)^{\lfloor(pn-1)/(p-1)\rfloor}}
\\\eq&\f{n_*!}{n_*+1}B_{p-1-n_*}\ (\mo\ p).
\endaligned\tag1.14$$
\endproclaim
\Proof. As $n-1$ is odd,  $pn-1\eq n-1\not\eq0\ (\mo\ p-1)$. Thus
$$\align
&p^n(-p)^{\lfloor(n-2)/(p-1)\rfloor}F_p(pn,0)
\\=&(-p)^{\lfloor(pn-1)/(p-1)\rfloor}F_p(pn,0)=C_p(pn,0)
\\=&\sum_{k=1}^n\bi{pn-1}{pk-1}(-1)^{pk}+\sum_{k=0}^{n-1}\bi {pn-1}{pk}(-1)^{pk}
\\=&\sum_{k=1}^n\bi{pn-1}{pk-1}(-1)^{k}+\sum_{k=0}^{n-1}\bi{pn-1}{p(n-k)-1}(-1)^{n-k}
\\=&2\sum_{k=1}^n\bi{pn-1}{pk-1}(-1)^{k}=2\sum_{k=1}^n\bi{pn-1}{p(n-k)}(-1)^{n-k}.
\endalign$$

Since
$$\align &p^{\lfloor(pn-1)/(p-1)\rfloor+1}=(1+(p-1))^{\lfloor(pn-1)/(p-1)\rfloor+1}
\\\gs&1+(p-1)\(\l\lfloor\f{pn-1}{p-1}\r\rfloor+1\)>1+(p-1)\f{pn-1}{p-1}=pn,
\endalign$$
we have
$$\l\lfloor\f{pn-2}{p-1}\r\rfloor=\l\lfloor\f{pn-1}{p-1}\r\rfloor
=n+\l\lfloor\f{n-1}{p-1}\r\rfloor\gs\ord_p(pn)$$
and hence
$$b-\l\lfloor\f{pn-1}{p-1}\r\rfloor>\l\lfloor\f{pn-1}{p-1}\r\rfloor\gs\ord_p(pn).$$
Recall that
$$(-1)^{pn-1}\f{(pn)!S(m\varphi(p^b),pn)}{(-p)^{\lfloor(pn-1)/(p-1)\rfloor}}
\eq F_p(pn,0)\ (\mo\ p^{b-\lfloor(pn-1)/(p-1)\rfloor})$$
and
$$\align&(-1)^{pn-2}\f{(pn-1)!S(m\varphi(p^b),pn-1)}{(-p)^{\lfloor(pn-2)/(p-1)\rfloor}}
\\\eq&F_p(pn-1,0)=\f{F_p(pn,0)}2\ (\mo\ p^{b-\lfloor(pn-1)/(p-1)\rfloor}).
\endalign$$

By the above, it suffices to show that
$$\f{F_p(pn,0)}{pn}\eq\f{n_*!}{n_*+1}B_{p-1-n_*}\ (\mo\ p).$$
Clearly $n^*\not=0,1$.
By Theorem 1.1 in the case $r=0$,
$$\align\f{F_p(p^2n,0)-F_p(pn,0)}{pn(p-1)}
\eq&\f{(-1)^{n^*}}{n^*!}\sum_{1<k\ls n^*}\bi{n^*}k\f{B_k}kB^{(0)}_{n^*-k}
\\\eq&-n_*!\f{B_{n^*}}{n^*}\eq\f{n_*!}{n_*+1}B_{p-1-n_*}\ (\mo\ p).\endalign$$
As $F_p(p^2n,0)/(pn)\eq0\ (\mo\ p)$ by Lemma 1.1, the desired result follows. \qed

\Remark\ 1.1. Let $p$ be a prime and let $n\in\Z^+$.

(i) If $p\not=2$ and $2\nmid n$, then $F_p(pn,0)=0$ because
$$C_p(pn,0)=\sum_{k=0}^{(n-1)/2}\(\bi{pn}{pk}(-1)^{pk}+\bi{pn}{p(n-k)}(-1)^{p(n-k)}\)=0.$$
If $m\in\Z^+$, $m\eq n\ (\mo\ p)$ and $m\eq n\not\eq0\ (\mo\ p-1)$, then by Theorem 1.2 and (1.4) we have
the following Kummer-type congruence:
$$\f{F_p(m,0)}m\eq \f{F_p(n,0)}n\ (\mo\ p).$$
If $1<n<p-1$ and $2\mid n$, then (1.13) yields (1.6).

(ii) When $2\mid n$ and $p-1\nmid n$, for any integers $m>0$ and $b>2\lfloor(pn-1)/(p-1)\rfloor$
we have
$$\ord_p\l((pn-1)!S(m\varphi(p^b),pn-1)\r)\gs\l\lfloor\f{pn-2}{p-1}\r\rfloor+\ord_p(pn)$$
and
$$\ord_p\l((pn)!S(m\varphi(p^b),pn)\r)\gs\l\lfloor\f{pn-1}{p-1}\r\rfloor+\ord_p(pn)$$
by (1.14). In 2001 I. M. Gessel and T. Lengyel [GL, Conjecture 1] conjectured that equality always holds
in our last two inequalities; this is not true since it might happen that
$B_{n^*}=B_{p-1-n_*}\eq0\ (\mo\ p)$, e.g., $\ord_{37}(B_{32})=\ord_{59}(B_{44})=1$.
($p$ is said to be {\it irregular} if $B_{2k}\eq0\ (\mo\ p)$
for some $0<k<(p-1)/2$. According to [IR, p.\,241] or [Mu, Theorem 2.13],
there are infinitely many irregular primes.)
By the way, Conjecture 2 of [GL] is an easy consequence of the congruence (1.4) due to Sun and Wan [SW2].

\proclaim{Corollary 1.4} Let $p\gs5$ be a prime.
Then
$$\bi{2p-1}{p-1}-1\eq-\f 23p^3B_{p-3}\ (\mo\ p^4)\tag1.15$$
and
$$\bi{4p}p-\bi{4p-1}{2p-1}-1\eq-\f{48}5p^5B_{p-5}\ (\mo\ p^6).\tag1.16$$
\endproclaim
\Proof. (1.15) follows from (1.13) in the case $n=2$.
Applying (1.13) with $n=4$ we find that
$$\f2{4p^5}\(-\bi{4p-1}{p-1}+\bi{4p-1}{2p-1}-\bi{4p-1}{3p-1}+\bi{4p-1}{4p-1}\)
\eq\f{4!}5B_{p-5}\ (\mo\ p).$$
This is equivalent to (1.16) since
$$\bi{4p-1}{p-1}+\bi{4p-1}{3p-1}=\bi{4p-1}{p-1}+\bi{4p-1}{p}=\bi{4p}p.$$
We are done. \qed

\Remark\ 1.2. (1.15) was first discovered by J.W.L. Glaisher (cf. [G1, p.\,21]
and [G2, p.\,323]). For a prime $p\gs5$ and $r\in\{1,\ldots,p-1\}$, we can determine
$\bi{2p-1}{p-1-r}$ mod $p^4$ in view of (1.15), and (1.1) and its proof.
\medskip

In the next theorem we determine $F_p(pn,r)$ mod $p^{\ord_p(n)+2}$
in the case $p-1\nmid n$ and $p\nmid r$.

\proclaim{Theorem 1.3} Let $p$ be an odd prime, and let $n>0$ and $r$
be integers with $p-1\nmid n$ and $p\nmid r$. Then, for any $b=1,\ldots,\ord_p(pn)$ we have
$$\aligned&\f{(-r)^nF_p(pn,r)+(-1)^{(b-1)n}\prod_{1\ls k\ls b'n_*,\,p\nmid k}k}{n_*!}
\\\eq&n_*(pB_{\varphi(p^b)}-p+1)-p^bn_*H_{n_*}+n(r^{p-1}-1)
\\&-pn\sum_{1<k<p-n_*}\bi{n_*+k}{n_*}\f{B_k}{kr^k}\ \ \ \ (\mo\ p^{b+1}),
\endaligned\tag1.17$$
where $b'=(p^b-1)/(p-1)$.
\endproclaim

\Remark\ 1.3. Let $p$ be an odd prime. A result of L. Carlitz [C]
states that $(B_k+p^{-1}-1)/k\in\Z_p$ for all $k\in\Z^+$ with $p-1\mid k$.
So we have $pB_{\varphi(p^b)}\eq p-1\ (\mo\ p^b)$ for all $b\in\Z^+$.

\proclaim{Corollary 1.5} Let $p$ be an odd prime. Let $n\in\Z^+$ and $r\in\Z$
with $p-1\nmid n$ and $p\nmid r$. Set $b=\ord_p(pn)$ and $b'=(p^b-1)/(p-1)$.
Then
$$F_p(pn,r)+\f{(-1)^{bn}}{r^n}\prod^{b'n_*}\Sb k=1\\p\nmid k\endSb k\in p^b\Z_p=pn\Z_p,\tag1.18$$
\endproclaim
\Proof. This is because the right-hand side of the congruence (1.17) belongs to $p^b\Z_p$. \qed

\proclaim{Corollary 1.6} Let $p$ be an odd prime. If $n\in\Z^+$, $r\in\Z$, $p-1\nmid n$ and $p\nmid r$,
then we can determine $F_p(pn,r)$ mod $p^2$ in the following way:
$$\aligned&\f{(-r)^nF_p(pn,r)+n_*!}{n_*!n_*}+pH_{n_*}
-pB_{p-1}+p-1
\\&\eq\f {pn}{n_*}\(q_p(r)-\sum_{1<k<p-n_*}\bi{n_*+k}{n_*}\f{B_k}{kr^k}\)\ (\mo\ p^2),
\endaligned\tag1.19$$
where $q_p(r)$ denotes the Fermat quotient $(r^{p-1}-1)/p$.
\endproclaim
\Proof. Just apply (1.17) with $b=1$. \qed

\Remark\ 1.4. If $p$ is a prime, $n\in\N$ and $r\in\Z$, then
$$F_p(pn+s,r)=\sum_{t=0}^s\bi st(-1)^tF_p(pn,r-t)\ \ \t{for any}\ s=0,\ldots,n^*;$$
because $\lfloor(pn+s-1)/(p-1)\rfloor=(pn+n^*)/(p-1)-1=\lfloor(pn-1)/(p-1)\rfloor$ and
$C_p(pn+s,r)$ coincides with
$$\sum_{k\eq r\,(\mo\ p)}\sum_{t=0}^s\bi st\bi{pn}{k-t}(-1)^k
=\sum_{t=0}^s\bi st(-1)^tC_p(pn,r-t)$$
by the Chu-Vandermonde convolution identity (cf. [GKP, (5.27)]).
\medskip

In the next section we determine $(\zeta_p^a-1)^{p^bn}$ modulo $p^{b+1}\pi^{p^bn}$ (where
$a\in\Z$ and $b\in\N$)
in terms of Bernoulli numbers or higher-order Bernoulli numbers.
On the basis of this, we prove Theorem 1.1 and Lemma 1.1 in Section 3 by a $p$-adic method.
In the proof of Theorem 1.3 given in Section 4, we have to employ the $p$-adic $\Gamma$-function
and the Gross-Koblitz formula for Gauss sums.
In Section 5, we study extended Fleck quotients and give an extension of (1.4)
which implies the following generalization of (1.5).

\proclaim{Theorem 1.4} Let $p$ be a prime, and let $a\in\Z^+$ and $l,m,n\in\N$ with
$m<p$ and $2\ls n-l-m\ls p$. Then we have
$$\aligned&\f1{p^{n-l}}\sum_{l<k\ls n}\bi{p^an-p^{a-1}m-d}{p^ak-p^{a-1}m-d}(-1)^{pk}\bi{k-1}l
\\&\qquad\eq\f{(-1)^{l-1}n!/l!}{\prod_{k=0}^m(n-l-k)}B^{(m+1)}_{p-n+l+m}\ (\mo\ p)
\endaligned\tag1.20$$
for all $d=1,\ldots,\max\{p^{a-2},1\}$.
\endproclaim

\heading{2. A theorem on roots of unity}\endheading

 In this section we establish the following auxiliary result.

\proclaim{Theorem 2.1} Let $p$ be a prime and define $\pi$ as in $(1.7)$.

{\rm (i)} We have
$$\pi^{p-1}\eq-p\ (\mo\ p^2),\ \ \t{i.e.,}\ \f{\pi^{p-1}}p\eq-1\ (\mo\ p).\tag2.1$$

{\rm (ii)} Let $a\in\Z$ and $n\in\N$.
If $m\in\N$ and $m\eq-n\ (\mo\ p)$, then
$$(\zeta_p^a-1)^n\eq\sum_{j=0}^{p-2}B_j^{(m)}\f{(a\pi)^{n+j}}{j!}\ (\mo\ p\pi^n).\tag2.2$$
For each $b\in\Z^+$ we have
$$\aligned(\zeta_p^a-1)^{p^bn}\eq&(a\pi)^{p^bn}+p^bn\sum_{1<k<p-1}\f{B_k}{k!k}(a\pi)^{p^bn+k}
\ (\mo\ p^{b+1}\pi^{p^bn}).
\endaligned\tag2.3$$
\endproclaim

\proclaim{Lemma 2.1} Let $p$ be any prime.
Then $\ord_p(\pi)=1/(p-1)$ and $\pi^{p-1}/p\eq-1\ (\mo\ \pi)$.
Also,
$$\zeta_p^a\eq\sum_{k=0}^{p-1}\f{(a\pi)^k}{k!}\ (\mo\ p\pi)\quad\t{for all}\ a\in\Z.\tag2.4$$
\endproclaim
\Proof. Clearly
$$\f{\pi}{\zeta_p-1}=\sum_{k=1}^{p-1}\f{(1-\zeta_p)^{k-1}}k=1-(1-\zeta_p)\eta$$
for some $\eta\in\Z_p[\zeta_p]$, hence
$$\f{\pi}{\zeta_p-1}\sum_{j=0}^{p-2}(1-\zeta_p)^j\eta^j=1-(1-\zeta_p)^{p-1}\eta^{p-1}.$$
Since $p/(1-\zeta_p)^{p-1}=\prod_{a=1}^{p-1}((1-\zeta_p^a)/(1-\zeta_p))$
is a unit in the ring $\Z_p[\zeta_p]$, by the above $\pi/(1-\zeta_p)$
and $\pi^{p-1}/p$ are also units in $\Z_p[\zeta_p]$ and hence
$\ord_p(\pi)=\ord_p(1-\zeta_p)=1/(p-1)$.
As $\pi/(1-\zeta_p)\eq -1\ (\mo\ 1-\zeta_p)$ and
$$\f p{(1-\zeta_p)^{p-1}}=\prod_{a=1}^{p-1}(1+\zeta_p+\cdots+\zeta_p^{a-1})
\eq\prod_{a=1}^{p-1}a\eq-1\ (\mo\ \zeta_p-1),$$
we have
$$\f{\pi^{p-1}}{p}=\f{(1-\zeta_p)^{p-1}}p\l(\f{\pi}{1-\zeta_p}\r)^{p-1}\eq-(-1)^{p-1}\eq-1\ (\mo\ \pi).$$

Write
$$\sum_{k=0}^{p-1}\f{(-\sum_{j=1}^{p-1}x^j/j)^k}{k!}=P(x)+x^pQ(x)$$
with $P(x),Q(x)\in\Z_p[x]$ and $\deg P(x)<p$.
If $-1<x<1$ then
$$1-x=e^{\log(1-x)}=\sum_{k=0}^\infty\f{(\log(1-x))^k}{k!}
=\sum_{k=0}^\infty\f{(-\sum_{j=1}^\infty{x^j}/j)^k}{k!}.$$
Comparing the coefficients of $1,x,\ldots,x^{p-1}$ we find that
$P(x)=1-x$. Therefore
$$\sum_{k=0}^{p-1}\f{\pi^k}{k!}=P(1-\zeta_p)+(1-\zeta_p)^pQ(1-\zeta_p)\eq\zeta_p\ (\mo\ \pi^p).$$

If $j\in\N$ and $\zeta_p^j\eq\sum_{k=0}^{p-1}(j\pi)^k/k!\ (\mo\ \pi^p)$, then
$$\align\zeta_p^{j+1}\eq&\sum_{k=0}^{p-1}\f{j^k\pi^k}{k!}\sum_{l=0}^{p-1}\f{\pi^l}{l!}
\\\eq&\sum_{n=0}^{p-1}\(\sum_{k=0}^n\bi nkj^k\)\f{\pi^n}{n!}=\sum_{n=0}^{p-1}\f{(j+1)^n\pi^n}{n!}\ (\mo\ \pi^p).
\endalign$$
Thus, by induction, $\zeta_p^a\eq\sum_{k=0}^{p-1}(a\pi)^k/k!\ (\mo\ \pi^p)$ for any $a\in\N$.

A general integer $a$ can be written in the form $pq+r$ with $q,r\in\Z$ and $0\ls r<p$.
In view of the above,
$$\align\zeta_p^a=\zeta_p^r\eq&\sum_{k=0}^{p-1}\f{(r\pi)^k}{k!}
=1+\sum_{k=1}^{p-1}\f{(a-pq)^k\pi^k}{k!}
\\\eq&1+\sum_{k=1}^{p-1}\f{a^k\pi^k}{k!}=\sum_{k=0}^{p-1}\f{(a\pi)^k}{k!}\ (\mo\ p\pi).
\endalign$$
This concludes the proof. \qed

\proclaim{Lemma 2.2} Let $k\in\Z^+$ and $m\in\N$. Then
$$\sum\Sb i_1,\ldots,i_k\in\N\\\sum_{j=1}^k i_jj=k\endSb
(-1)^{i_1+\cdots+i_k}\f{(\sum_{j=1}^k i_j-1)!}{i_1!\cdots i_k!}\prod_{j=1}^k\l(\f{B_j^{(m)}}{j!}\r)^{i_j}
=m\f{(-1)^{k}B_k}{k!k}.\tag2.5$$
\endproclaim
\Proof. For $0<x<2\pi$ we have
$$\align&\f{\t d}{\t dx}\(\log\f{e^x-1}x-\sum_{n=1}^\infty\f{B_n}n\cdot\f{(-x)^n}{n!}\)
\\=&\f{e^x}{e^x-1}-\f1x-\sum_{n=1}^\infty B_n\f{(-1)^nx^{n-1}}{n!}
\\=&\f1{1-e^{-x}}+\sum_{n=0}^\infty B_n\f{(-x)^{n-1}}{n!}
\\=&\f1{1-e^{-x}}+\f1{-x}\cdot\f{-x}{e^{-x}-1}=0.
\endalign$$
So $f(x)=\log((e^x-1)/x)-\sum_{n=1}^\infty (-x)^nB_n/(n!n)$ is a constant for $x\in(0,2\pi)$.
Letting $x\to0$ we find that the constant is zero.

 In light of the above,
 $$\align m\f{(-1)^{k-1}B_k}{k!k}=&-m[x^k]\log\l(\f{e^x-1}x\r)=[x^k]\log\l(\f x{e^x-1}\r)^{m}
 \\=&[x^k]\log\(1+\sum_{j=1}^\infty B^{(m)}_j\f{x^j}{j!}\)
 \\=&[x^k]\sum_{n=1}^k\f{(-1)^{n-1}}n\(\sum_{j=1}^k\f{B^{(m)}_j}{j!}x^j\)^n
 \\=&[x^k]\sum_{n=1}^k\f{(-1)^{n-1}}n\sum\Sb i_1,\ldots,i_k\in\N\\i_1+\cdots+i_k=n\endSb
 \f{n!}{i_1!\cdots i_k!}\prod_{j=1}^k\l(\f{B^{(m)}_j}{j!}x^j\r)^{i_j}
 \\=&\sum\Sb i_1,\ldots,i_k\in\N\\\sum_{j=1}^k i_jj=k\endSb
(-1)^{i_1+\cdots+i_k-1}\f{(\sum_{j=1}^k i_j-1)!}{i_1!\cdots i_k!}\prod_{j=1}^k\l(\f{B_j^{(m)}}{j!}\r)^{i_j}.
\endalign$$
So we have the desired (2.5). \qed
\medskip

\noindent{\it Proof of Theorem 2.1}.
(i) When $p=2$, (2.1) is trivial since $\pi=\zeta_2-1=-2$.

Now we consider the case $p>2$. In view of (2.4),
$$\align\sum_{a=1}^{p-1}\f{\zeta_p^a}{a^p}\eq&\sum_{k=0}^{p-1}\f{\pi^k}{k!}\sum_{a=1}^{p-1}\f1{a^{p-k}}
\\\eq&\sum_{a=1}^{(p-1)/2}\l(\f1{a^p}+\f1{(p-a)^p}\r)+\pi\sum_{a=1}^{p-1}\f1{a^{p-1}}
+\sum_{1<k<p}\f{\pi^k}{k!}\sum_{a=1}^{p-1}a^{k-1}
\\\eq&\sum_{a=1}^{(p-1)/2}\l(\f1{a^p}+\f1{(-a)^p}\r)+\pi(p-1)\eq-\pi\ (\mo\ p\pi).
\endalign$$
(It is well known that $\sum_{a=1}^{p-1}a^j\eq0\ (\mo\ p)$ for any $j\in\N$ with $p-1\nmid j$,
see, e.g., [IR, pp.\,235--236].)
For the norm
$$\al:=N_{\Q_p(\zeta_p)/\Q}\(\sum_{a=1}^{p-1}\f{\zeta_p^a}{a^p}\)\in\Z_p,$$
we have
$$\align \al=&\prod_{k=1}^{p-1}\sum_{a=1}^{p-1}\f{(\zeta_p^k)^a}{a^p}
=\prod_{k=1}^{p-1}\(k^p\sum_{a=1}^{p-1}\f{\zeta_p^{ka}}{(ka)^p}\)
\\\eq&\prod_{k=1}^{p-1}\(k^p\sum_{a=1}^{p-1}\f{\zeta_p^{ka}}{(\{ka\}_p)^p}\)
=((p-1)!)^p\prod_{k=1}^{p-1}\sum_{a=1}^{p-1}\f{\zeta_p^a}{a^p}
\\\eq&(-1)^p\(\sum_{a=1}^{p-1}\f{\zeta_p^a}{a^p}\)^{p-1}
\eq-\pi^{p-1}\(\f1{-\pi}\sum_{a=1}^{p-1}\f{\zeta_p^a}{a^p}\)^{p-1}\eq-\pi^{p-1}\ (\mo\ p^2).
\endalign$$
(Note that $(\beta+p\gamma)^p\eq\beta^p\ (\mo\ p^2)$ for any $\beta,\gamma\in\Z_p[\zeta_p]$.)
Thus
$\al\eq-\pi^{p-1}\eq p\ (\mo\ p\pi)$ and hence $\ord_p(\al-p)\gs\ord_p(p\pi)=1+1/(p-1)>1$.
As $\al-p\in\Z_p$, we must have $\ord_p(\al-p)\gs2$ and so $-\pi^{p-1}\eq\al\eq p\ (\mo\ p^2)$.
This proves (2.1).

(ii) Let $b,m\in\N$ and $m\eq-n\ (\mo\ p)$.
Observe that if $0<|x|<2\pi$ then
$$\align&\(\sum_{k=0}^{p-2}B_k^{(m)}\f{x^k}{k!}+\sum_{k=p-1}^{\infty}B_k^{(m)}\f{x^k}{k!}\)
\(\sum_{k=1}^{p-1}\f{x^{k-1}}{k!}+\sum_{k=p}^{\infty}\f{x^{k-1}}{k!}\)^{m}
\\&\qquad=\l(\f x{e^x-1}\r)^{m}\l(\f{e^x-1}x\r)^{m}=1.
\endalign$$
By comparing the coefficients of $1,x,\ldots,x^{p-2}$ we find that
$$\(\sum_{k=0}^{p-2}B^{(m)}_k\f{x^k}{k!}\)
\(\sum_{k=1}^{p-1}\f{x^{k-1}}{k!}\)^{m}=1+x^{p-1}Q(x)$$
for some $Q(x)\in\Z_p[x]$. It follows that
$$\(\sum_{k=0}^{p-2}B^{(m)}_k\f{(a\pi)^k}{k!}\)\(\sum_{k=1}^{p-1}\f{(a\pi)^{k-1}}{k!}\)^{m}\eq1\ (\mo\ p)$$
and
$$\(\sum_{k=0}^{p-2}B^{(m)}_k\f{(a\pi)^k}{k!}\)^{p^b}
\eq\(\sum_{k=1}^{p-1}\f{(a\pi)^{k-1}}{k!}\)^{-p^bm}\ (\mo\ p^{b+1}).$$
(Note that $(\beta+p^i\gamma)^p\eq\beta^p\ (\mo\ p^{i+1})$ for any $i\in\N$ and $\beta,\gamma\in\Z_p[\zeta_p]$.)

Since $\sum_{k=1}^{p-1}(a\pi)^{k-1}/k!=1+\pi\beta$ for some $\beta\in\Z_p[\zeta_p]$,
we have
$$\(\sum_{k=1}^{p-1}\f{(a\pi)^{k-1}}{k!}\)^{p^{b+1}}=(1+\pi\beta)^{p^{b+1}}\eq 1\ (\mo\ p^{b+1}\pi).$$
Note that $p^bn\eq -p^bm\ (\mo\ p^{b+1})$.
Therefore
$$\align\(\sum_{k=1}^{p-1}\f{(a\pi)^{k-1}}{k!}\)^{p^bn}\eq&\(\sum_{k=1}^{p-1}\f{(a\pi)^{k-1}}{k!}\)^{-p^bm}
\\\eq&\(\sum_{j=0}^{p-2}B^{(m)}_j\f{(a\pi)^j}{j!}\)^{p^b}
\ (\mo\ p^{b+1}).\endalign$$
In view of Lemma 2.1,
$$\f{\zeta_p^a-1}{\pi}\eq a\sum_{k=1}^{p-1}\f{(a\pi)^{k-1}}{k!}\ (\mo\ p).$$
Thus
$$\aligned\l(\f{\zeta_p^a-1}{\pi}\r)^{p^bn}\eq &\(a\sum_{k=1}^{p-1}\f{(a\pi)^{k-1}}{k!}\)^{p^bn}
\\\eq& a^{p^bn}\(\sum_{j=0}^{p-2}B^{(m)}_j\f{(a\pi)^j}{j!}\)^{p^b}\ (\mo\ p^{b+1}).
\endaligned\tag2.6$$
In the case $b=0$, this yields (2.2).

Below we assume  $b>0$ and want to prove (2.3).

By the multi-nomial theorem and the fact that $\pi^{2p-2}\eq0\ (\mo\ p^2)$,
$$\align&\(\sum_{j=0}^{p-2}B^{(m)}_j\f{(a\pi)^j}{j!}\)^p
\\=&\sum\Sb i_0,\ldots,i_{p-2}\in\N
\\i_0+\cdots+i_{p-2}=p\endSb\f {p!}{i_0!\cdots i_{p-2}!}\prod_{j=0}^{p-2}\l(B_j^{(m)}\f{(a\pi)^j}{j!}\r)^{i_j}
\\\eq&\sum_{k=0}^{2p-3}(a\pi)^k\sum\Sb \sum_{j=0}^{p-2}i_j=p\\\sum_{j=0}^{p-2}i_j j=k\endSb
\f {p!}{i_0!\cdots i_{p-2}!}\prod_{j=0}^{p-2}\l(\f{B_j^{(m)}}{j!}\r)^{i_j}\ (\mo\ p^2).
\endalign$$
If $i_0,\ldots,i_{p-2}\in\N$, $\sum_{j=0}^{p-2}i_j=p$
and $\sum_{j=0}^{p-2}i_j j=k\not\eq0\ (\mo\ p)$ then $i_0,\ldots,i_{p-2}$ are all smaller than $p$
and hence $p!/\prod_{j=0}^{p-2}i_j!\eq0\ (\mo\ p)$. Thus
$$\align&\(\sum_{j=0}^{p-2}B^{(m)}_j\f{(a\pi)^j}{j!}\)^p
\\\eq&(a\pi)^0\f{p!}{p!0!\cdots0!}\l(\f{B_0^{(m)}}{0!}\r)^p
\\&+[\![p\not=2]\!](a\pi)^p\f{p!}{0!p!0!\cdots 0!}\l(\f{B_1^{(m)}}{1!}\r)^p
\\&+\sum_{0<k<p-1}(a\pi)^k\sum\Sb \sum_{j=0}^{p-2}i_j=p\\\sum_{j=0}^{p-2}i_j j=k\endSb
\f {p!}{i_0!\cdots i_{p-2}!}\prod_{j=0}^{p-2}\l(\f{B_j^{(m)}}{j!}\r)^{i_j}
\ (\mo\ p^2).
\endalign$$
Note that $B_0^{(m)}=1$ and
$$B^{(m)}_1=[x]\l(\f x{e^x-1}\r)^{m}=[x]\(1-\f x2+\sum_{k=2}^\infty B_k\f{x^k}{k!}\)^{m}=-\f{m}2.$$
If $p\not=2$ then $(-m/2)^p\eq-m/2\eq n/2\ (\mo\ p)$.
If $0<k<p-1$, then
$$\align &\sum\Sb \sum_{j=0}^{p-2}i_j=p\\\sum_{j=0}^{p-2}i_j j=k\endSb
\f {p!}{i_0!\cdots i_{p-2}!}\prod_{j=0}^{p-2}\l(\f{B_j^{(m)}}{j!}\r)^{i_j}
\\=&\sum\Sb i_1,\ldots,i_{k}\in\N\\\sum_{j=1}^{k}i_j j=k\endSb\f{p!(B_0^{(m)}/0!)^{p-\sum_{j=1}^{k}i_j}}
{(p-\sum_{j=1}^{k}i_j)!}\prod_{j=1}^{k}\f{(B^{(m)}_j/j!)^{i_j}}{i_j!}
\\=&\sum\Sb i_1,\ldots,i_{k}\in\N\\\sum_{j=1}^{k}i_j j=k\endSb\ \prod_{0\ls i<\sum_{j=1}^{k}i_j}(p-i)
\times\prod_{j=1}^{k}\f{(B^{(p-n)}_j/j!)^{i_j}}{i_j!}
\\\eq&p\sum\Sb i_1,\ldots,i_{k}\in\N\\\sum_{j=1}^{k}i_j j=k\endSb(-1)^{\sum_{j=1}^{k}i_j-1}
\f{(\sum_{j=1}^{k}i_j-1)!}{i_1!\cdots i_k!}
\prod_{j=1}^{k}\l(\f{B^{(m)}_j}{j!}\r)^{i_j}
\ (\mo\ p^2).
\endalign$$
Therefore, with help of Lemma 2.2, we have
$$\(\sum_{j=0}^{p-2}B^{(m)}_j\f{(a\pi)^j}{j!}\)^p\eq1+p\pi S\ (\mo\ p^2)$$
where
$$\align S=&[\![p\not=2]\!]\f n2\cdot\f{a^p\pi^{p-1}}{p}-m\sum_{0<k<p-1}a^k\pi^{k-1}\f{(-1)^kB_k}{k!k}
\\\eq&[\![p\not=2]\!]\f{an}2\l(\f{\pi^{p-1}}p+1\r)+n\sum_{1<k<p-1}a^k\pi^{k-1}\f{B_k}{k!k}
\\\eq&n\sum_{1<k<p-1}a^k\pi^{k-1}\f{B_k}{k!k}\ \ (\mo\ p)\quad\t{(by (2.1))}.
\endalign$$
(Recall that $B_1=-1/2$ and $B_{2j+1}=0$ for all $j\in\Z^+$.)

Since $b>0$, it follows from the above that
$$\align\(\sum_{j=0}^{p-2}B_j^{(m)}\f{(a\pi)^j}{j!}\)^{p^b}
\eq&(1+p\pi S)^{p^{b-1}}
\\\eq&1+p^{b-1}p\pi S+\sum_{1<k\ls p^{b-1}}p^{b-1}\bi{p^{b-1}-1}{k-1}\f{p^k\pi^k}kS^k
\\\eq&1+p^b\pi S\  (\mo\ p^{b+1}),
\endalign$$
where we have noted that $\pi^2/2\in\Z_p[\zeta_p]$ and $p^{k-2}/k\in\Z_p$ for $k=3,4,\ldots$
(cf. [S03, Lemma 2.1]). Combining this with (2.6),
we find that
$$\align (\zeta_p^a-1)^{p^bn}\eq&(a\pi)^{p^bn}(1+p^b\pi S)
\\\eq&(a\pi)^{p^bn}+p^bn\sum_{1<k<p-1}(a\pi)^{p^bn+k}\f{B_k}{k!k}\ (\mo\ p^{b+1}\pi^{p^bn}).
\endalign$$
This proves (2.3) and we are done.  \qed

\heading{3. Proofs of Theorem 1.1 and Lemma 1.1}\endheading

\proclaim{Lemma 3.1} Let $p$ be a prime, and let $n\in\N$ and $r\in\Z$. Then
$$pC_p(n,r)=\sum_{a=0}^{p-1}\zeta_p^{-ar}(1-\zeta_p^a)^n.\tag3.1$$
\endproclaim
\Proof. This known result can be easily proved
by using $[\![p\mid k-r]\!]=p^{-1}\sum_{a=0}^{p-1}\zeta_p^{a(k-r)}$.
\qed

\medskip
\noindent
{\it Proof of Lemma 1.1}. For $a=1,\ldots,p-1$,
as $(\zeta_p^a-1)/\pi\eq a\ (\mo\ \pi)$ by Lemma 2.1,
we have
$$\l(\f{\zeta_p^a-1}{\pi}\r)^n\eq a^n\ \l(\mo\ p^{\ord_p(n)}\pi\r)$$
since
$$\l(\f{\zeta_p^a-1}{\pi}\r)^p\eq a^p\ (\mo\ p\pi),\
\l(\f{\zeta_p^a-1}{\pi}\r)^{p^2}\eq a^{p^2}\ (\mo\ p^2\pi),\ \ldots.$$
Let $g$ be a primitive root modulo $p$. Then $g^n\not\eq1\ (\mo\ p)$ (as $p-1\nmid n$),
and also
$$\align(g^n-1)\sum_{a=1}^{p-1}a^n=&\sum_{a=1}^{p-1}(ag)^n-\sum_{a=1}^{p-1}a^n
\\\eq&\sum_{a=1}^{p-1}(\{ag\}_p)^n-\sum_{a=1}^{p-1}a^n=0\ (\mo\ p^{\ord_p(n)+1}).
\endalign$$
Therefore $\sum_{a=1}^{p-1}a^n\eq0\ (\mo\ p^{\ord_p(n)+1})$
and hence
$$\sum_{a=1}^{p-1}\l(\f{\zeta_p^a-1}{\pi}\r)^n
\eq\sum_{a=1}^{p-1}a^n\eq0\ \l(\mo\ p^{\ord_p(n)}\pi\r).$$

On the other hand,
$$\sum_{a=0}^{p-1}\l(\f{\zeta_p^a-1}{\pi}\r)^n=\f {pC_p(n,0)}{(-\pi)^n}$$
by Lemma 3.1.
So we have
$$\ord_p(C_p(n,0))\gs\ord_p\l(p^{\ord_p(n)-1}\pi^{n+1}\r)
=\ord_p(n)-1+\f{n+1}{p-1}$$
and hence
$$\ord_p(F_p(n,0))=\ord_p(C_p(n,0))-\l\lfloor\f{n-1}{p-1}\r\rfloor>\ord_p(n)-1.$$
Since $F_p(n,0)\in\Z$, this shows that $F_p(n,0)/n\in\Z_p$. We are done. \qed

\proclaim{Lemma 3.2} Let $p$ be a prime, and let $n\in\N$, $r\in\Z$ and $r\not\eq0\ (\mo\ p)$. Then
$$\sum_{a=1}^{p-1}a^n(\zeta_p^{ar}-1)\eq-\f{(r\pi)^{n^*}}{n^*!}p^{[\![p-1\mid n]\!]}
\eq n_*!(-r\pi)^{n^*}p^{[\![p-1\mid n]\!]}\ (\mo\ p\pi).\tag3.2$$
\endproclaim
\Proof. In view of Lemma 2.1,
$$\sum_{a=1}^{p-1}a^n(\zeta_p^{ar}-1)\eq\sum_{a=1}^{p-1}a^n\sum_{k=1}^{p-1}\f{(ar\pi)^k}{k!}
=\sum_{k=1}^{p-1}\f{(r\pi)^k}{k!}\sum_{a=1}^{p-1}a^{n+k}\ (\mo\ p\pi).$$
Since $\sum_{a=1}^{p-1}a^{n+k}\eq-[\![p-1\mid n+k]\!]\ (\mo\ p)$, we have
$$\align\sum_{a=1}^{p-1}a^n(\zeta_p^{ar}-1)\eq&-\sum\Sb 1\ls k\ls p-1\\p-1\mid k-n^*\endSb\f{(r\pi)^k}{k!}
=-\f{(r\pi)^{n^*}}{n^*!}\l(\f{(r\pi)^{p-1}}{(p-1)!}\r)^{[\![n^*=0]\!]}
\\\eq&-\f{(r\pi)^{n^*}}{n^*!}(-\pi^{p-1})^{[\![n^*=0]\!]}
\eq-\f{(r\pi)^{n^*}}{n^*!}p^{[\![p-1\mid n]\!]}
\\\eq&n_*!(-r\pi)^{n^*}p^{[\![p-1\mid n]\!]}\ \ (\mo\ p\pi).
\endalign$$
This proves (3.2). \qed

\medskip
\noindent{\it Proof of Theorem 1.1}. Without loss of generality,
we assume $m>n$ and write $m-n=p^b(p-1)d$ with $b,d\in\Z^+$ and $p\nmid d$.
Clearly $2\mid p^b(p-1)$.
Set $$D=\f1{\pi^{n^*}}\sum_{a=1}^{p-1}\zeta_p^{-ar}\l(\f{\zeta_p^a-1}{\pi}\r)^n
\(\l(\f{\zeta_p^a-1}{\pi}\r)^{p^b(p-1)d}-1\).$$
Then
$$\align(-1)^nD=&\f{(-1)^{n+p^b(p-1)d}}{\pi^{p^b(p-1)d+n+n^*}}
\sum_{a=1}^{p-1}\zeta_p^{-ar}(\zeta_p^a-1)^{n+p^b(p-1)d}
\\&-\f{(-1)^n}{\pi^{n+n^*}}\sum_{a=1}^{p-1}\zeta_p^{-ar}(\zeta_p^a-1)^n
\\=&\f{pC_p(n+p^b(p-1)d,r)}{\pi^{p^b(p-1)d+n+n^*}}-\f{pC_p(n,r)-[\![n=0]\!]}{\pi^{n+n^*}}
\ \ (\t{by Lemma 3.1})
\\=&\f{(-p)^{(n+n^*)/(p-1)}}{\pi^{n+n^*}}F_p(n,r)
-\f{(-p)^{p^bd+(n+n^*)/(p-1)}}{\pi^{p^b(p-1)d+n+n^*}}F_p(m,r)
\endalign$$
and hence
$$\align(-1)^nD\l(\f{\pi^{p-1}}{-p}\r)^{(n+n^*)/(p-1)}
=&F_p(n,r)-\l(\f{-p}{\pi^{p-1}}\r)^{p^bd}F_p(m,r)
\\\eq& F_p(n,r)-F_p(m,r)\ (\mo\ p^b\pi).
\endalign$$
(Note that $(-p/\pi^{p-1})^{p^b}\eq1\ (\mo\ p^b\pi)$ since $-p/\pi^{p-1}\eq1\ (\mo\ \pi)$.)

Let $a$ be an integer not divisible by $p$.
In view of Theorem 2.1(ii),
$$\l(\f{\zeta_p^a-1}{\pi}\r)^n\eq\sum_{j=0}^{p-2}a^{n+j}B_j^{(\{-n\}_p)}\f{\pi^j}{j!}\ (\mo\ p),$$
and
$$\align&\l(\f{\zeta_p^a-1}{\pi}\r)^{p^b(p-1)d}-1
\\\eq&a^{p^b(p-1)d}-1+p^b(p-1)d\sum_{1<k<p-1}\f{B_k}{k!k}a^{p^b(p-1)d+k}\pi^k
\\\eq&-p^bd\sum_{1<k<p-1}\f{B_k}{k!k}a^k\pi^k\ (\mo\ p^{b+1})
\endalign$$
with help of the congruence
$a^{\varphi(p^{b+1})}\eq1\ (\mo\ p^{b+1})$ (Euler's theorem).

In light of the above,
$$ \pi^{n^*}D\eq-p^bd\sum_{j=0}^{p-2}
\sum_{1<k<p-1}\f{B_j^{(\{-n\}_p)}B_k}{j!k!k}\pi^{j+k}G_r(n+j+k)\ (\mo\ p^{b+1}),$$
where $G_r(s)=\sum_{a=1}^{p-1}a^s\zeta_p^{-ar}$ for $s\in\Z$.
By Lemma 3.2, for $k=0,\ldots,p-2$ we have
$$\align G_r(n+k)=&\sum_{a=1}^{p-1}a^{n+k}(\zeta_p^{-ar}-1)+\sum_{a=1}^{p-1}a^{n+k}
\\\eq&[\![k<n^*]\!](n_*+k)!(r\pi)^{n^*-k}-[\![k=n^*]\!]
\\\eq&[\![k\ls n^*]\!](n_*+k)!(r\pi)^{n^*-k}\ (\mo\ \pi^{n^*+1}).
\endalign$$
(Note that $(n_*+n^*)!=(p-1)!\eq-1\ (\mo\ p)$ by Wilson's theorem.)
So $\pi^{n^*}D$ is congruent to
$$-p^bd\sum_{1<k\ls n^*}\sum_{j=0}^{n^*-k}\f{B_j^{(\{-n\}_p)}B_k}{j!k!k}
\pi^{j+k}(n_*+j+k)!(r\pi)^{n^*-j-k}$$
modulo $p^b\pi^{n^*+1}$ and hence
$$\align D\eq&-p^bd\sum_{1<l\ls n^*}(n_*+l)!r^{n^*-l}
\sum_{1<k\ls l}\f{B_k}{k!k}\cdot\f{B^{(\{-n\}_p)}_{l-k}}{(l-k)!}
\\\eq&-p^bd\times n_*!\Sigma \ \ \ (\mo\ p^b\pi),
\endalign$$
where $$\Sigma=\sum_{1<l\ls n^*}\bi{n_*+l}{n_*}r^{n^*-l}
\sum_{1<k\ls l}\bi lk\f{B_k}kB^{(\{-n\}_p)}_{l-k}.$$
Therefore
$$F_p(m,r)-F_p(n,r)\eq(-1)^{n-1}D
\eq p^bd(-1)^nn_*!\Sigma\ (\mo\ p^b\pi),$$
i.e., the $p$-adic order of the rational number
$$R=\f{F_p(m,r)-F_p(n,r))}{m-n}+(-1)^{n_*}n_*!\Sigma$$
is at least $\ord_p(\pi)=1/(p-1)>0$. It follows that $\ord_p(R)\gs 1$.

If $0<l\ls n^*$, then
$$\align\bi{n_*+l}{n_*}=&\prod_{k=1}^l\f{p-1-n^*+k}{k}
\\\eq&(-1)^l\prod_{k=1}^l\f{n^*-k+1}{j}=(-1)^l\bi{n^*}l\ (\mo\ p).
\endalign$$
Note also that
$$\align&\sum_{1<l\ls n^*}(-1)^l\bi{n^*}lr^{n^*-l}
\sum_{1<k\ls l}\bi lk\f{B_k}kB^{(\{-n\}_p)}_{l-k}
\\=&(-1)^{n^*}\sum_{1<k\ls n^*}\bi{n^*}k\f{B_k}k
\sum_{l=k}^{n^*}\bi{n^*-k}{l-k}B^{(\{-n\}_p)}_{l-k}(-r)^{n^*-k-(l-k)}
\\=&(-1)^{n^*}\sum_{1<k\ls n^*}\bi{n^*}k\f{B_k}kB^{(\{-n\}_p)}_{n^*-k}(-r).
\endalign$$
So we have
$$\align\f{F_p(m,r)-F_p(n,r))}{m-n}\eq&(-1)^{n_*-1}n_*!\Sigma\eq\f{\Sigma}{n^*!}
\\\eq&\f{(-1)^{n^*}}{n^*!}\sum_{1<k\ls n^*}\bi{n^*}k\f{B_k}kB^{(\{-n\}_p)}_{n^*-k}(-r)\ (\mo\ p).
\endalign$$
This concludes the proof. \qed

\heading{4. Proof of Theorem 1.3}\endheading

The following lemma in the case $b=1$ is a known result due to
Beeger in 1913 (cf. [Mu, p.\,23]).

\proclaim{Lemma 4.1} Let $p$ be an odd prime, and let $b$ be a positive integer.
Then
$$w_{p^b}\eq \f{pB_{\varphi(p^b)}-p+1}{p^b}\ (\mo\ p),\tag4.1$$
where $w_{p^b}$ denotes the generalized Wilson quotient $(1+\prod_{0<a<p^b,\,p\nmid a}a)/p^b$.
\endproclaim
\Proof. Let $g$ be a primitive root modulo $p^b$. As Gauss discovered,
 $$\prod^{p^b-1}\Sb a=1\\p\nmid a\endSb a\eq\prod_{k=0}^{\varphi(p^b)-1}g^k=g^{\varphi(p^b)(\varphi(p^b)-1)/2}
 \eq g^{\varphi(p^b)/2}\eq-1\ (\mo\ p^b).$$
 So $w_{p^b}$ is an integer.

 Clearly
 $$\f{3B_{\varphi(3)}-3+1}3=\f{3/6-3+1}3=-\f12\eq w_3=\f{2!+1}3\ (\mo\ 3).$$
 Below we assume $p^b>3$.

Let $k>1$ be an integer. Recall that
$$\align \sum_{a=1}^{p-1}a^k=&\f{B_{k+1}(p)-B_{k+1}}{k+1}
=\f1{k+1}\sum_{j=0}^k\bi{k+1}{j+1}B_{k-j}p^{j+1}
\\=&pB_k+pk\sum_{j=1}^k\bi{k-1}{j-1}(pB_{k-j})\f{p^{j-1}}{j(j+1)}.
\endalign$$
Since $p\not=2$, by [S03, Lemma 2.1] we have $p^{j-2}/(j(j+1))\in\Z_p$ for $j=3,4,\ldots$.
Note also that $pB_{k-j}\in\Z_p\ (0\ls j\ls k)$
by the von Staudt-Clausen theorem.
So
$$\sum_{a=1}^{p-1}a^k\eq pB_k+pk\l(\f p2B_{k-1}+(k-1)pB_{k-2}\f p{2\times 3}\r)
\ (\mo\ p^{\ord_p(k)+2}).$$
When $p-1\mid k$, we have $B_{k-1}\in\Z_p$ and
$pB_{k-2}\eq-[\![p=3\ \&\ k\not=2]\!]\ (\mo\ p)$
by the von Staudt-Clausen theorem, therefore
$$\sum_{a=1}^{p-1}a^k\eq pB_k-[\![p=3\ \&\ k\not=2]\!]p\f{k(k-1)}2\ (\mo\ p^{\ord_p(k)+2}).$$
Putting $k=\varphi(p^b)>2$, we obtain
$$\align\sum_{a=1}^{p-1}a^{\varphi(p^b)}\eq& pB_{\varphi(p^b)}-[\![p=3]\!]p\varphi(p^b)\f{\varphi(p^b)-1}2
\\\eq& pB_{\varphi(p^b)}+[\![p=3]\!]p^b\ (\mo\ p^{b+1}).\endalign$$
(Note that if $p=3$ then $b>1$ and hence $p\mid\varphi(p^b)$.)
Thus
$$\align\sum_{a=1}^{p-1}\f{a^{\varphi(p^b)}-1}{p^b}\eq&\f{pB_{\varphi(p^b)}+[\![p=3]\!]p^b-p+1}{p^b}
\\\eq&\f{pB_{\varphi(p^b)}-p+1}{p^b}+[\![p=3]\!]\ (\mo\ p).
\endalign$$

If $a_1$ and $a_2$ are two integers relatively prime to $p^b$, then
$$\align\f{(a_1a_2)^{\varphi(p^b)}-1}{p^b}
=&\f{a_1^{\varphi(p^b)}-1}{p^b}a_2^{\varphi(p^b)}+\f{a_2^{\varphi(p^b)}-1}{p^b}
\\\eq&\f{a_1^{\varphi(p^b)}-1}{p^b}+\f{a_2^{\varphi(p^b)}-1}{p^b}\ (\mo\ p^b)
\endalign$$
by Euler's theorem. Therefore
$$\align \sum^{p^b-1}\Sb a=1\\p\nmid a\endSb\f{a^{\varphi(p^b)}-1}{p^b}
\eq&\f{(\prod_{0<a<p^b,\,p\nmid a}a)^{\varphi(p^b)}-1}{p^b}=\f{(-1+p^bw_{p^b})^{\varphi(p^b)}-1}{p^b}
\\\eq&\f{(1-\varphi(p^b)p^bw_{p^b})-1}{p^b}=-\varphi(p^b)w_{p^b}\eq p^{b-1}w_{p^b}\ (\mo\ p^b).
\endalign$$

Suppose that $n$ is an integer with $0<n<b$. Then
$$\align&\sum^{p^{n+1}-1}\Sb a=1\\p\nmid a\endSb\f{a^{\varphi(p^b)}-1}{p^b}
=\sum^{p^n-1}\Sb a=1\\p\nmid a\endSb\sum_{s=0}^{p-1}\f{(p^ns+a)^{\varphi(p^b)}-1}{p^b}
\\=&\f1{p^b}\sum^{p^n-1}\Sb a=1\\p\nmid a\endSb\sum_{s=0}^{p-1}\(a^{\varphi(p^b)}-1
+\sum_{k=1}^{\varphi(p^b)}\f{\varphi(p^b)}k\bi{\varphi(p^b)-1}{k-1}(p^ns)^ka^{\varphi(p^b)-k}\)
\\=&p\sum^{p^n-1}\Sb a=1\\p\nmid a\endSb\f{a^{\varphi(p^b)}-1}{p^b}
+(p-1)\sum_{k=1}^{\varphi(p^b)}\bi{\varphi(p^b)-1}{k-1}\f{p^{kn-1}}k
\sum_{s=1}^{p-1}s^k\sum^{p^n-1}\Sb a=1\\p\nmid a\endSb a^{\varphi(p^b)-k}.
\endalign$$
If $n>1$ and $k\gs2$, then
$$\f{p^{kn-1}}k=p^{k(n-1)+1}\f{p^{k-2}}k\eq0\ (\mo\ p^{n+1})$$
because $k(n-1)\gs2(n-1)\gs n$ and $p^{k-2}/k\in\Z_p$ by [S03, Lemma 2.1].
In the case $n=1$, as $p^{k-3}/k\in\Z_p$ for $k=4,5,\ldots$ (cf. [S03, Lemma 2.1]) and
$$\f{p^2}3\sum_{s=1}^{p-1}s^3=\f{p^2}3\sum_{s=1}^{(p-1)/2}(s^3+(p-s)^3)\eq0\ (\mo\ p^2),$$
we have
$$\f{p^{k-1}}k\sum_{s=1}^{p-1}s^k\eq0\ (\mo\ p^2)\quad\t{for}\ k=3,4,5,\ldots.$$
Thus
$$\align&\sum^{p^{n+1}-1}\Sb a=1\\p\nmid a\endSb\f{a^{\varphi(p^b)}-1}{p^b}
-p\sum^{p^n-1}\Sb a=1\\p\nmid a\endSb\f{a^{\varphi(p^b)}-1}{p^b}
\\\eq&(p-1)p^{n-1}\sum_{s=1}^{p-1}s\sum^{p^n-1}\Sb a=1\\p\nmid a\endSb a^{\varphi(p^b)-1}
\\&+[\![n=1]\!](p-1)(\varphi(p^b)-1)\f p2\sum_{s=1}^{p-1}s^2\sum_{a=1}^{p-1}a^{\varphi(p^b)-2}
\\\eq&(p-1)p^{n-1}\f{p(p-1)}2\sum^{(p^n-1)/2}\Sb a=1\\p\nmid a\endSb\l(\f1a+\f1{p^n-a}\r)
\\&+[\![n=1]\!](p-1)(\varphi(p^b)-1)\f p2\cdot\f{p(p-1)(2p-1)}6\sum_{a=1}^{p-1}a^{\varphi(p^b)-2}
\\\eq&[\![p=3\ \&\ n=1]\!]\f{p}2\times\f{2\cdot3\cdot5}6\sum_{a=1}^21\eq-p[\![p=3\ \&\ n=1]\!]
\ (\mo\ p^{n+1}).
\endalign$$

In view of the above,
$$\align &p^{b-1}\sum_{a=1}^{p-1}\f{a^{\varphi(p^b)}-1}{p^b}
-\sum^{p^b-1}\Sb a=1\\p\nmid a\endSb\f{a^{\varphi(p^b)}-1}{p^b}
\\=&\sum_{0<n<b}p^{b-n-1}\(p\sum^{p^n-1}\Sb a=1\\p\nmid a\endSb\f{a^{\varphi(p^b)}-1}{p^b}
-\sum^{p^{n+1}-1}\Sb a=1\\p\nmid a\endSb\f{a^{\varphi(p^b)}-1}{p^b}\)
\\\eq&\sum_{0<n<b}p^{b-n-1}p[\![p=3\ \&\ n=1]\!]=p^{b-1}[\![p=3]\!]\ (\mo\ p^b)
\endalign$$
and hence
$$p^{b-1}w_{p^b}\eq\sum^{p^b-1}\Sb a=1\\p\nmid a\endSb\f{a^{\varphi(p^b)}-1}{p^b}
\eq p^{b-1}\sum_{a=1}^{p-1}\f{a^{\varphi(p^b)}-1}{p^b}-p^{b-1}[\![p=3]\!]\ (\mo\ p^b).$$
Therefore
$$w_{p^b}\eq\sum_{a=1}^{p-1}\f{a^{\varphi(p^b)}-1}{p^b}-[\![p=3]\!]\eq\f{pB_{\varphi(p^b)}-p+1}{p^b}\ (\mo\ p).$$
This concludes the proof. \qed

\proclaim{Lemma 4.2} Let $p$ be an odd prime and let $n\in\Z^+$, $r\in\Z$ and $r\not\eq0\ (\mo\ p)$.
Then, for any $b=1,\ldots,\ord_p(pn)$, we have
$$\aligned&\f{r^n}{(-p)^N}\sum_{a=1}^{p-1}(a\pi)^{pn}\zeta_p^{-ar}
-(-1)^{(b-1)n}\prod^{b'n_*}\Sb k=1\\p\nmid k\endSb k
\\\eq&n_*!\l(p^b n_*H_{n_*}-n_*(pB_{\varphi(p^b)}-p+1)-pnq_p(r)\r)
\ (\mo\ p^b\pi),\endaligned\tag4.2$$
where
$$b'=\f{p^b-1}{p-1}\ \ \t{and}\ \ N=\f{pn+n^*}{p-1}=\l\lfloor\f{pn-1}{p-1}\r\rfloor+1.$$
\endproclaim
\Proof. Write $pn=p^bm$ with $m\in\Z^+$. Then
$$\align& r^{(p-1)n}=(1+pq_p(r))^{p^{b-1}m}
\\=&1+p^bmq_p(r)+\sum_{1<k\ls p^{b-1}m}p^{b-1}m\bi{p^{b-1}m-1}{k-1}\f {p^k}kq_p(r)^k
\\\eq&1+p^bmq_p(r)\ \ (\mo\ p^{b+1})
\endalign$$
since $p^{k-2}/k\in\Z_p$ for $k=2,3,\ldots$. Thus
$$\align\sum_{a=1}^{p-1}a^{pn}\zeta_p^{-ar}=&\f{(-1)^{pn}r^{-n}}{r^{(p-1)n}}
\sum_{a=1}^{p-1}(-ar)^{p^bm}\zeta_p^{-ar}
\eq\f{(-1)^{n}r^{-n}}{1+p^bmq_p(r)}
\sum_{s=1}^{p-1}s^{p^bm}\zeta_p^s
\\\eq&\f{(-1)^n}{r^n}(1-p^bmq_p(r))\sum_{a=1}^{p-1}a^{pn}\zeta_p^a\ (\mo\ p^{b+1}).
\endalign$$

 Let $\omega$ be the Teichm\"uler character
of the multiplicative group
$$(\Z/p\Z)^*=\{\bar a=a+p\Z:\,a=1,\ldots,p-1\}.$$
Then for each $a=1,\ldots,p-1$ the value
$\omega(\bar a)$ is just the unique $(p-1)$-th root of unity
(in the algebraic closure of $\Q_p$) with $\omega(\bar a)\eq a\ (\mo\ p)$.
(See, e.g., [Wa, p.\,51].)
Since $a^{p^b}\eq\omega(\bar a)^{p^b}\ (\mo\ p^{b+1})$, we have
$$\sum_{a=1}^{p-1}a^{pn}\zeta_p^a\eq\sum_{a=1}^{p-1}\omega(\bar a)^{pn}\zeta_p^a
=\sum_{a=1}^{p-1}\omega(\bar a)^{-n^*}\zeta_p^a\ (\mo\ p^{b+1}).$$
By the Gross-Koblitz formula for Gauss sums (cf. [BEW, p.\,350] and [GK]),
$$G(n^*):=\sum_{a=1}^{p-1}\omega(\bar a)^{-n^*}\zeta_p^a
=-\pi_0^{n^*}\Gamma_p\l(\f{n^*}{p-1}\r)$$
where $\Gamma_p$ is Morita's $p$-adic $\Gamma$-function (see [BEW, p.\,277]
or [Mu, p.\,59 and pp. 67--70] for the definition and basic properties), and
$\pi_0$ is the unique element in $\Z_p[\zeta_p]$ satisfying
$$\pi_0^{p-1}=-p\ \ \t{and}\ \ \pi_0\eq\zeta_p-1\ (\mo\ (\zeta_p-1)^2).$$
(See [Go, pp.\,172--173] for the existence of $\pi_0$.)
Clearly $\pi_0\eq\zeta_p-1\eq\pi\ (\mo\ \pi^2)$ and hence $\pi_0/\pi\eq1\ (\mo\ \pi)$.
(Furthermore, we have $\pi_0/\pi=(\pi_0/\pi)^p\pi^{p-1}/(-p)\eq1\ (\mo\ p)$.)

In view of the above,
$$\sum_{a=1}^{p-1}a^{pn}\zeta_p^{-ar}
\eq\f{(-1)^{n-1}}{r^n}\pi_0^{n^*}(1-pnq_p(r))\Gamma_p\l(\f{n^*}{p-1}\r)\ (\mo\ p^{b+1})$$
and hence
$$\align &(-1)^{n-1}r^n\sum_{a=1}^{p-1}(a\pi)^{pn}\zeta_p^{-ar}
\\\eq&\pi_0^{pn+n^*}\l(\f{\pi}{\pi_0}\r)^{pn}
(1-pnq_p(r))\Gamma_p\l(\f{n^*}{p-1}\r)\ (\mo\ p^{b+1}\pi^{pn})
\\\eq&(-p)^N\l(\f{\pi}{\pi_0}\r)^{p^bm}
(1-pnq_p(r))\Gamma_p\l(\f{n^*}{p-1}\r)\ (\mo\ p^b\pi^{pn+n^*+1})
\\\eq&(-p)^N(1-pnq_p(r))\Gamma_p\l(\f{n^*}{p-1}\r)\ (\mo\ p^{b+N}\pi).
\endalign$$
(Note that $(\pi/\pi_0)^{p^b}\eq1\ (\mo\ p^b\pi)$.)

Since
$$\f{n^*}{p-1}=1-\f{n_*}{p-1}\eq1+n_*\f{p^{b+1}-1}{p-1}=1+n_*+pn_*+\cdots+p^bn_*\ (\mo\ p^{b+1}),$$
we have $$\Gamma_p\l(\f{n^*}{p-1}\r)\eq\Gamma_p((p^b+b')n_*+1)
=(-1)^{(p^b+b')n_*+1}\prod^{(p^b+b')n_*}\Sb k=1\\p\nmid k\endSb k\ (\mo\ p^{b+1}).$$
Observe that
$$\align \prod^{(p^b+b')n_*}\Sb k=1\\p\nmid k\endSb k
=&\prod_{s=0}^{n_*-1}\prod^{p^b-1}\Sb t=1\\p\nmid t\endSb (p^bs+t)
\times\prod^{b'n_*}\Sb k=1\\p\nmid k\endSb(p^bn_*+k)
\\=&\(\prod^{p^b-1}\Sb t=1\\p\nmid t\endSb t\)^{n_*}
\prod_{s=0}^{n_*-1}\prod^{p^b-1}\Sb t=1\\p\nmid t\endSb\l(1+p^b\f st\r)
\times\prod^{b'n_*}\Sb k=1\\p\nmid k\endSb k\times\prod^{b'n_*}\Sb k=1\\p\nmid k\endSb\l(1+p^b\f{n_*}k\r)
\endalign$$
and hence
$$\align &\(\prod^{(p^b+b')n_*}\Sb k=1\\p\nmid k\endSb k\)\bigg/
\(\(\prod^{p^b-1}\Sb t=1\\p\nmid t\endSb t\)^{n_*}\prod^{b'n_*}\Sb k=1\\p\nmid k\endSb k\)
\\\eq&1+p^b\sum_{s=0}^{n_*-1}\sum^{p^b-1}\Sb t=1\\p\nmid t\endSb\f st+p^bn_*\(\sum_{0<k<(b'-1)n_*}\f1k
+\sum_{j=1}^{n_*}\f1{n_*(b'-1)+j}\)
\\\eq&1+p^bn_*\sum_{j=1}^{n_*}\f1{j}=1+p^bn_*H_{n_*}\ (\mo\ p^{b+1})
\endalign$$
since $p\mid b'-1$ and $\sum_{k=1}^{p-1}1/k=\sum_{0<k<p/2}(1/k+1/(p-k))\eq0\ (\mo\ p)$.

By Lemma 4.1,
 $$\align\(-\prod^{p^b-1}\Sb t=1\\p\nmid t\endSb t\)^{n_*}=&(1-p^bw_{p^b})^{n_*}\eq1-n_*p^bw_{p^b}
 \\\eq&1-n_*(pB_{\varphi(p^b)}-p+1)\ (\mo\ p^{b+1}).
 \endalign$$
Therefore
 $$\align &\Gamma_p\l(\f {n_*}{p-1}\r)\bigg/\prod^{b'n_*}\Sb k=1\\p\nmid k\endSb k
 \\\eq&(-1)^{(p^b+b')n_*+1}(1+p^bn_*H_{n_*})(-1)^{n_*}(1-n_*(pB_{\varphi(p^b)}-p+1))
 \\\eq&(-1)^{bn+1}(1+p^bn_*H_{n^*}-n_*(pB_{\varphi(p^b)}-p+1))\ (\mo\ p^{b+1}).
 \endalign$$

 Note that $b'-1=p\sum_{0\ls i<b-1}p^i$ and hence
 $$\align \prod^{b'n_*}\Sb k=1\\p\nmid k\endSb k=&\prod_{0\ls s<n_*(b'-1)/p}\prod_{t=1}^{p-1}(ps+t)
 \times\prod_{j=1}^{n_*}(n_*(b'-1)+j)
 \\\eq&((p-1)!)^{n_*(b'-1)/p}n_*!\eq(-1)^{n(b-1)}n_*!\ (\mo\ p).
 \endalign$$
 So there is a $u\in\Z$ such that
 $$U:=(-1)^{(b-1)n}\prod^{b'n_*}\Sb k=1\\p\nmid k\endSb k=n_*!+pu.$$
 Combining the above, we finally get
 $$\align&\f{r^n}{(-p)^N}\sum_{a=1}^{p-1}(a\pi)^{pn}\zeta_p^{-ar}
 \\\eq&(-1)^{n-1}(1-pnq_p(r))\Gamma_p\l(\f{n_*}{p-1}\r)
 \\\eq&(-1)^{n-1}(1-p^bmq_p(r))\times(-1)^{bn+1}\prod^{b'n_*}\Sb k=1\\p\nmid k\endSb k
 \\&\times\(1+p^bn_*\l(H_{n_*}-\f{pB_{\varphi(p^b)}-p+1}{p^b}\r)\)
 \\\eq&(n_*!+pu)\(1-p^bmq_p(r)+p^bn_*\l(H_{n_*}-\f{pB_{\varphi(p^b)}-p+1}{p^b}\r)\)
 \\\eq&U+n_*!\l(-pnq_p(r)+p^bn_*H_{n^*}-n_*(pB_{\varphi(p^b)}-p+1)\r)
 \ (\mo\ p^{b}\pi).
 \endalign$$
 This yields the desired (4.2). \qed

\medskip
\noindent{\it Proof of Theorem 1.3}. Let $\bar n=pn$. By Lemma 3.1 and Theorem 2.1(ii),
$$\align (-1)^{p\bar n}pC_p(p\bar n,r)=&\sum_{a=0}^{p-1}\zeta_p^{-ar}(\zeta_p^a-1)^{p\bar n}
\\\eq&\sum_{a=1}^{p-1}\zeta_p^{-ar}(a\pi)^{p\bar n}\ \l(\mo\ p^{\ord_p(p\bar n)}\pi^{p\bar n}\r).
\endalign$$
As $b\ls\ord_p(\bar n)$ and $p\bar n-pn=(p-1)\bar n\eq0\ (\mo\ \varphi(p^{b+1}))$, we have
$$(-1)^npC_p(p\bar n,r)\eq\pi^{(p-1)\bar n}\sum_{a=1}^{p-1}(a\pi)^{pn}\zeta_p^{-ar}
\ (\mo\ p^{b+1}\pi^{p\bar n}).$$
Set $N=(\bar n+n^*)/(p-1)$. Then
$$N+\bar n=\f{p\bar n+n^*}{p-1}=\l\lfloor\f{p\bar n-1}{p-1}\r\rfloor+1.$$
Thus
$$\f{(-1)^n}{(-p)^N}\sum_{a=1}^{p-1}(a\pi)^{pn}\zeta_p^{-ar}
=\f{pC_p(p\bar n,r)}{(-p)^{N+\bar n}}\l(\f{-p}{\pi^{p-1}}\r)^{\bar n}
\eq-F_p(p\bar n,r)\ (\mo\ p^b\pi).$$
(Note that $(-p/\pi^{p-1})^{\bar n}\eq 1\ (\mo\ p^b\pi)$ since
$-p/\pi^{p-1}\eq1\ (\mo\ \pi)$ and $p^b\mid \bar n$.)
Combining this with Lemma 4.2, we obtain
$$\align&-(-r)^nF_p(p\bar n,r)-(-1)^{(b-1)n}\prod^{b'n_*}\Sb k=1\\p\nmid k\endSb k
\\\eq&n_*!\l(p^bn_*H_{n_*}-n_*(pB_{\varphi(p^b)}-p+1)-pnq_p(r)\r)\ (\mo\ p^b\pi)
\endalign$$
and hence
$$\aligned&\f{(-r)^nF_p(p\bar n,r)+(-1)^{(b-1)n}\prod_{1\ls k\ls b'n_*,\,p\nmid k}k}{n_*!}
\\\eq&n_*(pB_{\varphi(p^b)}-p+1)-p^bn_*H_{n_*}+pnq_p(r)\ (\mo\ p^{b+1}).
\endaligned\tag4.3$$

By Theorem 1.1,
$$\f{F_p(p\bar n,r)-F_p(\bar n,r)}{(p-1)\bar n}
\eq(-1)^{n_*-1}\f{n_*!}{r^n}\sum_{1<k\ls n^*}\bi{n_*+k}{n_*}\f{B_k}{kr^k}\ (\mo\ p)$$
and so
$$\f{(-r)^n(F_p(p\bar n,r)-F_p(\bar n,r))}{n_*!}
\eq\bar n\sum_{1<k\ls n^*}\bi{n_*+k}{n_*}\f{B_k}{kr^k}\ (\mo\ p^{b+1}).\tag4.4$$
From (4.3) and (4.4) we immediately get the desired congruence (1.17).
\qed

\heading{5. Congruences for extended Fleck quotients}\endheading

Let $p$ be a prime, and let $a\in\Z^+$, $n\in\N$ and $r\in\Z$.
In 1977 C. S. Weisman [We] extended Fleck's inequality
by showing that
$$\ord_p\(\sum_{k\eq r\,(\mo\ p^a)}\bi nk(-1)^k\)\gs\l\lfloor\f{n-p^{a-1}}{\varphi(p^a)}\r\rfloor.$$
An extension of this result was given by the author in [S06].
During his study of the $\psi$-operator in Fontaine's theory in 2005,
D. Wan finally obtained the following extension of Fleck's inequality (cf. [W] and [SW1]):
For any $l\in\N$ we have
$$\ord_p\(\sum_{k\eq r\,(\mo\ p^a)}\bi nk(-1)^k\bi{(k-r)/p^a}l\)
\gs\l\lfloor\f{n-lp^a-p^{a-1}}{\varphi(p^a)}\r\rfloor,$$
i.e., the extended Fleck quotient
$$F^{(l)}_{p^a}(n,r):=(-p)^{\lfloor(n-lp^a-p^{a-1})/\varphi(p^a)\rfloor}
\sum_{k\eq r\,(\mo\ p^a)}\bi nk(-1)^k\bi{(k-r)/p^a}l$$
is an integer. In this section we study $F^{(l)}_{p^a}(n,r)$ mod $p$.

\proclaim{Theorem 5.1} Let $p$ be a prime, and let $a\in\Z^+$, $l,n\in\N$, $r\in\Z$
and $s,t\in\{0,\ldots,p^{a-1}-1\}$.
Let $m\in\N$ with $m\eq-n\ (\mo\ p)$. Then
$$\aligned&(-1)^{l+t-1}F^{(l)}_{p^a}(p^{a-1}n+s,p^{a-1}r+t)
\\\eq&[\![n>l]\!]\bi st\bi{\lfloor(n-l-1)/(p-1)\rfloor}l(n-l)_*B^{(m)}_{(n-l)^*}(-r)\ (\mo\ p),
\endaligned\tag5.1$$
provided that we have one of the
following {\rm (i)--(iii):}

\ \ {\rm (i)} $a=1\ \t{or}\ p\mid n\ \t{or}\ p-1\nmid n-l-1$;

\ \ {\rm (ii)} $\lfloor s/p^{a-2}\rfloor=2\lfloor t/p^{a-2}\rfloor$ and $p\not=2$;

\ \ {\rm (iii)} $\lfloor s/p^{a-2}\rfloor=\lfloor t/p^{a-2}\rfloor=p-1$.

\endproclaim

Now we deduce Theorem 1.4 from Theorem 5.1.

\medskip
\noindent{\it Proof of Theorem 1.4}.
Let $d\in\Z^+$ with $d\ls\max\{p^{a-2},1\}$. Then
$$\l\lfloor\f{(p^an-p^{a-1}m-d)-lp^a-p^{a-1}}{\varphi(p^a)}\r\rfloor
=\l\lfloor\f{p(n-l)-m-2}{p-1}\r\rfloor=n-l.$$
Also, $\lfloor(p^{a-1}-d)/p^{a-2}\rfloor=p-1$ if $a>1$.
Thus
$$\align &\f1{(-p)^{n-l}}\sum_{l<k\ls n}\bi{p^an-p^{a-1}m-d}{p^ak-p^{a-1}m-d}(-1)^{p^ak-p^{a-1}m-d}\bi{k-1}l
\\=&F^{(l)}_{p^a}(p^an-p^{a-1}m-d,p^a-p^{a-1}m-d)
\\=&F^{(l)}_{p^a}\l(p^{a-1}(pn-m-1)+p^{a-1}-d,p^{a-1}(p-m-1)+p^{a-1}-d\r)
\\\eq&(-1)^{l+(p^{a-1}-d)-1}\bi{p^{a-1}-d}{p^{a-1}-d}\bi{\lfloor(pn-m-1-l-1)/(p-1)\rfloor}l
\\&\times(pn-m-1-l)_*!B^{(m+1)}_{(pn-m-1-l)^*}(m+1-p)\quad\t{(by Theorem 5.1)}
\\\eq&(-1)^{l+d}\bi nl(n-m-1-l)!B^{(m+1)}_{p-1-(n-m-1-l)}(m+1-p)\ (\mo\ p).
\endalign$$
Since $B^{(m+1)}_0,\ldots,B^{(m+1)}_{p-2}\in\Z_p$ and
$$\align&(-1)^{p-n+l+m}B^{(m+1)}_{p-n+l+m}(m+1-p)=B^{(m+1)}_{p-n+l+m}(p)
\\=&\sum_{j=0}^{p-n+l+m}\bi{p-n+l+m}jB^{(m+1)}_jp^{p-n+l+m-j}\eq B^{(m+1)}_{p-n+l+m}\ (\mo\ p),
\endalign$$
by the above we have
$$\align&\f1{(-p)^{n-l}}\sum_{l<k\ls n}\bi{p^an-p^{a-1}m-d}{p^ak-p^{a-1}m-d}(-1)^{p^ak-p^{a-1}m-d}\bi{k-1}l
\\\eq&(-1)^{l+d}\f{n!/l!}{\prod_{k=0}^m(n-l-k)}(-1)^{p-n+l+m}B^{(m+1)}_{p-n+l+m}\ (\mo\ p),
\endalign$$
which is equivalent to (1.20). \qed
\smallskip

To prove Theorem 5.1 we need some lemmas.

\proclaim{Lemma 5.1} Let
$f(x)$ be a function from $\Z$ to a field, and let $m,n\in\Z^+$.
Then, for any $r\in\Z$ we have
$$\sum_{k=0}^n\bi nk(-1)^kf\l(\l\lfloor\f{k-r}m\r\rfloor\r)
=\sum_{k\eq \bar r\,(\mo\ m)}\bi{n-1}k(-1)^{k-1}\Delta f\l(\f{k-\bar r}m\r),$$
where $\bar r=r+m-1$ and $\Delta f(x)=f(x+1)-f(x)$.
\endproclaim
\Proof. This is Lemma 2.1 of Sun [S06]. \qed

\proclaim{Lemma 5.2} Let $p$ be a prime, and let $l,n\in\N$ with
$n>p$. Then
$$\aligned &F^{(l)}_p(n,r)+[\![l>0]\!]F^{(l-1)}_p(n-p,r)
\\\eq&-\sum_{k=1}^{p-1}\f1k\sum_{j=0}^{k-1}F^{(l)}_p(n-p+1,r-j)\ (\mo\ p).
\endaligned\tag5.2$$
\endproclaim
\Proof. Set $n'=n-(p-1)>0$. With help of the Chu-Vandermonde convolution identity,
$$\align &F^{(l)}_p(n,r)
\\=&(-p)^{-\lfloor(n-lp-1)/(p-1)\rfloor}\sum_{k\eq r\,(\mo\ p)}
\sum_{j=0}^{p-1}\bi{p-1}j\bi{n'}{k-j}(-1)^k\bi{(k-r)/p}l
\\=&-\f1p\sum_{j=0}^{p-1}\bi{p-1}j
(-p)^{-\lfloor(n'-lp-1)/(p-1)\rfloor}\sum_{p\mid k-r}\bi{n'}{k-j}(-1)^k\bi{(k-r)/p}l
\\=&-\f1p\sum_{j=0}^{p-1}\bi{p-1}j(-1)^jF^{(l)}_p(n',r-j).
\endalign$$

For any $j=0,\ldots,p-1$,
clearly
$$\align&\bi{p-1}j(-1)^j=\prod_{0<i\ls j}\l(1-\f pi\r)
\\\eq&1-\sum_{0<i\ls j}\f pi\eq(-1)^{p-1}+p\sum_{j<k<p}\f 1k\ (\mo\ p^2).
\endalign$$
(Recall that $H_{p-1}=\sum_{k=1}^{p-1}1/k\eq0\ (\mo\ p)$ if $p\not=2$.)
Also,
$$\align&-\f1p\sum_{j=0}^{p-1}F^{(l)}_p(n',r-j)
\\=&(-p)^{-1-\lfloor(n'-lp-1)/(p-1)\rfloor}
\sum_{k=0}^{n'}\bi{n'}k(-1)^k\bi{\lfloor(k-r+p-1)/p\rfloor}l
\\=&(-p)^{-\lfloor((n'-1)-(l-1)p-1)/(p-1)\rfloor}\sum_{k\eq r\,(\mo\ p)}\bi{n'-1}k(-1)^{k-1}\bi{(k-r)/p}{l-1}
\\=&-[\![l>0]\!]F^{(l-1)}_p(n'-1,r),
\endalign$$
where we have applied Lemma 5.1 with $f(x)=\bi xl$ for the second equality and view $\bi x{-1}$ as $0$.
Therefore
$$\align F^{(l)}_p(n,r)\eq&(-1)^p[\![l>0]\!]F^{(l-1)}_p(n'-1,r)-\sum_{j=0}^{p-1}\sum_{j<k<p}\f{F^{(l)}_p(n',r-j)}k
\\\eq&-[\![l>0]\!]F^{(l-1)}_p(n'-1,r)-\sum_{k=1}^{p-1}\f1k\sum_{j=0}^{k-1}F^{(l)}_p(n',r-j)\ (\mo\ p).
\endalign$$
This proves (5.2). \qed

\proclaim{Lemma 5.3} Let $p$ be a prime, and let $l,n\in\N$ and $r\in\Z$. If $n>lp$, then
$$F^{(l)}_p(n,r)\eq(-1)^l\bi{\lfloor(n-l-1)/(p-1)\rfloor}lF_p(n-lp,r)\ (\mo\ p).\tag5.3$$
\endproclaim
\Proof. We use induction on $l+n$.

Clearly $l=0$ and $n=1$ if $l+n=1$.
In the case $l=0$, (5.3) holds trivially for $n>0$.

Below we let $l>0$ and assume the corresponding result for smaller values of $l+n$.
As $n>lp$, we have $n'-1>(l-1)p$ where $n'=n-p+1$. By the induction hypothesis,
$(-1)^{l-1}F^{(l-1)}_p(n'-1,r)$ is congruent to
$$\align &\bi{\lfloor (n'-1-(l-1)-1)/(p-1)\rfloor}{l-1}F_p(n'-1-(l-1)p,r)
\\&\quad=\bi{\lfloor (n'-l-1)/(p-1)\rfloor}{l-1}F_p(n-lp,r)\endalign$$
modulo $p$.

Clearly $n'>lp-p+1\gs l$. If $n'\ls lp$ then
$$\f{n'-l-1}{p-1}-l=\f{n'-lp-1}{p-1}<0$$
and $$\sum_{k=1}^{p-1}\f1k\sum_{j=0}^{k-1}(-1)^lF^{(l)}_p(n',r-j)\eq0\ (\mo\ p).$$
If $n'>lp$, then by the induction hypothesis,
$$\align&\sum_{k=1}^{p-1}\f1k\sum_{j=0}^{k-1}(-1)^lF^{(l)}_p(n',r-j)
\\\eq&\sum_{k=1}^{p-1}\f1k\sum_{j=0}^{k-1}\bi{\lfloor (n'-l-1)/(p-1)}l F_p(n'-lp,r-j)
\\\eq&-\bi{\lfloor (n'-l-1)/(p-1)\rfloor}l F_p(n-lp,r)\ (\mo\ p),
\endalign$$
where we have applied Lemma 5.2 for Fleck quotients.

The above, together with Lemma 5.2, yields that
$$\align (-1)^l F^{(l)}_p(n,r)\eq&(-1)^{l-1}[\![l>0]\!]F^{(l-1)}_p(n'-1,r)
\\&-\sum_{k=1}^{p-1}\f1k\sum_{j=0}^{k-1}(-1)^{l}F^{(l)}_p(n',r-j)
\\\eq&\bi{\lfloor (n'-l-1)/(p-1)\rfloor}{l-1}F_p(n-lp,r)
\\&+\bi{\lfloor (n'-l-1)/(p-1)\rfloor}{l}F_p(n-lp,r)
\\\eq&\bi{\lfloor (n-l-1)/(p-1)\rfloor}{l}F_p(n-lp,r)\ (\mo\ p).
\endalign$$

The induction proof is now complete. \qed

\proclaim{Lemma 5.4} Let $p$ be a prime, and let $a\in\Z^+$, $n\in\N$, $r\in\Z$
and $s,t\in\{0,\ldots,p^{a-1}-1\}$. If one of {\rm (i)-(iii)} in Theorem 5.1 is satisfied, then
$$F^{(l)}_{p^a}(p^{a-1}n+s,p^{a-1}r+t)\eq(-1)^t\bi stF^{(l)}_p(n,r)\ (\mo\ p)\tag5.4$$
\endproclaim
\Proof. (5.4) holds trivially in the case $a=1$. Below we assume $a\gs2$.

Write $s=\sum_{k=0}^{a-2}s_kp^k$ and $t=\sum_{k=0}^{a-2}t_kp^k$ with $s_k,t_k\in\{0,\ldots,p-1\}$.
By [SW1, Theorem 1.1], if $a>2$ then
$$\align &F^{(l)}_{p^a}(p^{a-1}n+s,p^{a-1}r+t)
\\=&F^{(l)}_{p^a}\(p\(p^{a-2}n+\sum_{k=1}^{a-2}s_kp^{k-1}\)+s_0,
p\(p^{a-2}r+\sum_{k=1}^{a-2}t_kp^{k-1}\)+t_0\)
\\\eq&(-1)^{t_0}\bi{s_0}{t_0}F^{(l)}_{p^{a-1}}
\(p^{a-2}n+\sum_{k=1}^{a-2}s_kp^{k-1},p^{a-2}r+\sum_{k=1}^{a-2}t_kp^{k-1}\)
\\\eq&\cdots\eq\(\prod_{k=0}^{a-3}(-1)^{t_k}\bi{s_k}{t_k}\)F^{(l)}_{p^2}(pn+s_{a-2},pr+t_{a-2})\ (\mo\ p).
\endalign$$

Observe that $s_{a-2}=\lfloor s/p^{a-2}\rfloor$ and $t_{a-2}=\lfloor t/p^{a-2}\rfloor$.
If (i) or (ii) holds, then
$$F^{(l)}_{p^2}(pn+s_{a-2},pr+t_{a-2})\eq(-1)^{t_{a-2}}\bi{s_{a-2}}{t_{a-2}}F^{(l)}_p(n,r)\ (\mo\ p)\tag5.5$$
by [SW1, Theorem 1.2]. Suppose that (iii) holds (i.e., $s_{a-2}=t_{a-2}=p-1$) but (i) fails.
By [SW1, Lemma 3.3],
$$\align&(-1)^{\lfloor(pn+s_{a-2}-(n-1)p-1)/(p-1)\rfloor}
F^{(n-1)}_p(pn+s_{a-2},t_{a-2})
\\&\qquad\eq(-1)^{n+t_{a-2}}n\bi{s_{a-2}}{t_{a-2}}\f{\sigma}p\ \ (\mo\ p),
\endalign$$
where
$$\align \sigma=&1+(-1)^p\f{\prod_{i=2}^p(p(n-1)+p-1+i)}{\prod_{i=1}^{p-1}i}
=1+(-1)^p\prod_{k=1}^{p-1}\l(1+\f{pn}k\r)
\\\eq&1+(-1)^p\(1+pn\sum_{k=1}^{p-1}\f1k\)\eq0\ \ (\mo\ p^2).
\endalign$$
(Note that if $p=2$ then $n$ is odd since (i) fails.)
Thus $F^{(n-1)}_p(pn+s_{a-2},t_{a-2})\eq0\ (\mo\ p)$ and hence (5.5) holds by [SW1, Lemma 3.2].

Provided (i) or (ii) or (iii), by the above we have
$$\align &F^{(l)}_{p^a}(p^{a-1}n+s,p^{a-1}r+t)
\\\eq&\prod_{k=0}^{a-2}(-1)^{t_k}\bi{s_k}{t_k}\times F^{(l)}_p(n,r)
\\\eq&(-1)^{\sum_{k=0}^{a-2}t_kp^k}
\bi{\sum_{k=0}^{a-2}s_kp^k}{\sum_{k=0}^{a-2}t_kp^k}F^{(l)}_p(n,r)=(-1)^t\bi stF^{(l)}_p(n,r)\ (\mo\ p),
\endalign$$
where we have applied Lucas' theorem (cf. [HS]).
This completes the proof. \qed

\medskip
\noindent{\it Proof of Theorem 5.1}. In view of Lemma 5.4, it suffices to show that
$(-1)^lF^{(l)}_p(n,r)$ is congruent to
$$-[\![n>l]\!]\bi{\lfloor(n-l-1)/(p-1)\rfloor}l(n-l)_*B^{(m)}_{(n-l)^*}(-r)$$
modulo $p$.
In the case $n\ls lp$, this is easy since the last expression vanishes.

 Below we assume $n>lp$. By Lemma 5.3 and (1.4),
$$\align &(-1)^lF^{(l)}_p(n,r)
\\\eq&\bi{\lfloor(n-l-1)/(p-1)\rfloor}lF_p(n-lp,r)
\\\eq&-\bi{\lfloor(n-l-1)/(p-1)\rfloor}l(n-lp)_*!B^{(m)}_{(n-lp)^*}(-r)\ (\mo\ p).
\endalign$$
Since $(n-lp)^*=(n-l)^*$ and $(n-lp)_*=(n-l)_*$, the desired result follows and we are done. \qed

\bigskip

\Ack. This paper is based on the previous work [SW2] joint with
D. Wan. The author is indebted to Prof. Wan for his work in [SW2],
and Prof. K. Ono for his comments on the author's related talk
given at the University of Wisconsin at Madison in April 2006.
The author also thanks his two Ph.D. students H. Pan and H. Q. Cao for
discussion on a particular case of (1.20).

\enddocument